\documentclass[a4paper,oneside,10pt]{article}%
\usepackage{CJK,amsmath}
\usepackage{amsfonts}
\usepackage{amsthm}
\usepackage{indentfirst}
\usepackage{slashbox}
\usepackage{makeidx}
\usepackage{mathrsfs}
\usepackage{amsmath}
\usepackage{amssymb}
\usepackage{amsfonts}
\usepackage{undertilde}
\usepackage[perpage,symbol*]{footmisc}
\usepackage{bm}
\usepackage{graphicx}
\usepackage[T1]{fontenc}
\usepackage{titlesec}
\usepackage{ccaption}
\usepackage[tight,hang]{subfigure}
\usepackage{graphicx}%
\providecommand{\U}[1]{\protect\rule{.1in}{.1in}}

\input{undertilde}
\makeatletter

\newcommand{\Rmnum}[1]{\expandafter\@slowromancap\romannumeral #1@}
\makeatother
\newcommand*{\HEI}
{\CJKfamily{hei}}

\newtheorem{theorem}{Theorem}[section]
\newtheorem{lemma}{Lemma}[section]

\newtheorem{corollary}{Corollary}[section]
\newtheorem{definition}{Definition}[section]

\newtheorem{assumption}{Assumption}[section]
\begin{document}

\title{Non-Markovian Fully Coupled Forward-Backward Stochastic Systems and Classical
Solutions of Path-dependent PDEs\thanks{This work was supported by National
Natural Science Foundation of China (No. 11171187, No. 10871118 and No.
10921101).}}
\author{Shaolin Ji\thanks{Institute for Financial Studies and Institute of
Mathematics, Shandong University, Jinan, Shandong 250100, PR China
(Jsl@sdu.edu.cn, Fax: +86 0531 88564100).}
\and Shuzhen Yang\thanks{School of mathematics, Shandong University, Jinan, Shandong
250100, PR China. (yangsz@mail.sdu.edu.cn). }}
\date{}
\maketitle

\textbf{Abstract}. This paper explores the relationship between non-Markovian
fully coupled forward-backward stochastic systems and {path-dependent PDEs.
The definition of classical solution for the path-dependent PDE is given
within the framework of functional It\^{o} calculus. Under mild hypotheses, we
prove that the }forward-backward stochastic system provides the unique
classical solution to the path-dependent PDE.

\bigskip

{\textbf{Keywords}: Functional It\^{o} calculus, Non-Markovian fully coupled
forword-backward systems, Path-dependent PDEs, Classical solutions}

\addcontentsline{toc}{section}{\hspace*{1.8em}Abstract}

\section{Introduction}

Linear Backward Stochastic Differential Equations (BSDEs) was introduced by
Bismut \cite{Bismut}. The existence and uniqueness theorem for nonlinear BSDEs
was established by Pardoux and Peng \cite{Pardoux.E 3}. Then Peng \cite{Peng S
5} and Pardoux and Peng \cite{Pardoux.E 2} gave a relationship between
Markovian {forword-backward systems and }systems of quasilinear parabolic
PDEs, which generalized the classical Feynman-Kac formula. Peng \cite{Peng S
1} pointed out that for non-Markovian {forword-backward systems, }it was an
open problem to find the corresponding "PDE"{.}

Recently in the framework of {functional It\^{o} calculus, }a path-dependent
PDE was introduced by Dupire \cite{Dupire.B} which shed light on this problem
(for a recent account of this theory we refer the reader to \cite{Cont.R},
\cite{Cont-2} and \cite{Cont-3}). Inspired by Dupire's work, Peng and Wang
\cite{Peng S 3} obtained a nonlinear Feynman-Kac formula for classical
solutions of path-dependent PDEs in terms of non-Markovian BSDEs. Cosso
\cite{Cosso} proved that, under some assumptions, the non-Markovian
{forword-backward system gives the unique continuous viscosity solution to the
path-dependent PDE. For the further development, the readers may refer to
\cite{Peng S 2} and \cite{E}.}

In this paper, we study the relationship between solutions of {non-Markovian
fully coupled forword-backward systems and classical solutions of
}path-dependent PDEs. More precisely, the {non-Markovian forword-backward
system is described by the following fully coupled forword-backward} SDE:
\begin{equation}
X^{\gamma_{t},x}(s)=x+\int_{t}^{s}b(W_{s}^{\gamma_{t}},X^{\gamma_{t}%
,x}(s),Y^{\gamma_{t},x}(s),Z^{\gamma_{t},x}(s))ds+\int_{t}^{s}\sigma
(W_{s}^{\gamma_{t}},X^{\gamma_{t},x}(s),Y^{\gamma_{t},x}(s),Z^{\gamma_{t}%
,x}(s))dW(s), \tag{1.1}%
\end{equation}%
\begin{equation}
Y^{\gamma_{t},x}(s)=g(W_{T}^{\gamma_{t}},X^{\gamma_{t},x}(T))-\int_{s}%
^{T}h(W_{s}^{\gamma_{t}},X^{\gamma_{t},x}(s),Y^{\gamma_{t},x}(s),Z^{\gamma
_{t},x}(s))ds-\int_{s}^{T}Z^{\gamma_{t},x}(s)dW(s),\quad s\in\lbrack t,T].
\tag{1.2}%
\end{equation}
We first give {the definition of classical solution, within the framework of
functional It\^{o} calculus, for the path-dependent PDE. }Note that in our
context, we use the terminology "classical solution" to distinguish it from
"viscosity solution" in \cite{Peng S 2}, \cite{E} and \cite{Cosso}. {Then
}under mild hypotheses, we establish some estimates and regularity results for
the solution of the above system with respect to paths. Finally, we show that
the solution of (1.2) is related to the classical solution of the following
path-dependent PDE%
\[%
\begin{array}
[c]{l}%
D_{t}u(\gamma_{t},x)+\mathcal{L}u(\gamma_{t},x)+tr[\nabla_{x}D_{z}u(\gamma
_{t},x)\sigma(\gamma_{t},x,u(\gamma_{t},x),v(\gamma_{t},x))]+\frac{1}%
{2}tr[D_{zz}u(\gamma_{t},x)]\\
=h(\gamma_{t},x,u(\gamma_{t},x),v(\gamma_{t},x)),\\
v(\gamma_{t},x)=\nabla_{x}u(\gamma_{t},x)\sigma(\gamma_{t},x,u(\gamma
_{t},x),v(\gamma_{t},x))+D_{z}u(\gamma_{t},x),\\
u(\gamma_{T},x)=g(\gamma_{T},x),\quad\gamma_{T}\in{\Lambda}^{n},
\end{array}
\]
where
\[
\mathcal{L}u=\frac{1}{2}tr[(\sigma\sigma^{T})(\gamma_{t},x,u,v)\nabla
_{xx}u]+\langle b(\gamma_{t},x,u,v)\nabla_{x}u\rangle.
\]

The paper is organized as follows: in section 2, we present some fundamental
results on functional It\^{o} calculus and FBSDE theory. We establish some
estimates and regularity results for the solution of non-Markovian FBSDEs in
section 3. Finally, in section 4, we give the relationship between
non-Markovian fully coupled FBSDEs and {path-dependent PDEs.}

\section{Preliminaries}

\subsection{Functional It\^{o} calculus}

The following notations and tools are mainly from Dupire \cite{Dupire.B}. Let
$T>0$ be fixed. For each $t\in\lbrack0,T]$, we denote by $\Lambda_{t}$ the set
of c\`{a}dl\`{a}g $\mathbb{R}^{d}$-valued functions on $[0,t]$. For each
$\gamma\in\Lambda_{T}$ the value of $\gamma$ at time $s\in\lbrack0,T]$ is
denoted by $\gamma(s)$. Thus $\gamma=\gamma(s)_{0\leq s\leq T}$ is a
c\`{a}dl\`{a}g process on $[0,T]$ and its value at time $s$ is $\gamma(s)$.
The path of $\gamma$ up to time $t$ is denoted by $\gamma_{t}$, i.e.,
$\gamma_{t}=\gamma(s)_{0\leq s\leq t}\in\Lambda_{t}$. We denote $\Lambda
=\bigcup_{t\in\lbrack0,T]}\Lambda_{t}$. For each $\gamma_{t}\in\Lambda$ and
$x\in\mathbb{R}^{d}$ we denote by $\gamma_{t}(s)$ the value of $\gamma_{t}$ at
$s\in\lbrack0,t]$ and $\gamma_{t}^{x}:=(\gamma_{t}(s)_{0\leq s<t},\gamma
_{t}(t)+x)$ which is also an element in $\Lambda_{t}$.

Let $|\cdot|$ denote the norm in $\mathbb{R}^{d}$. We now define a distance on
$\Lambda$. For each $0\leq t,\bar{t}\leq T$ and $\gamma_{t},\bar{\gamma}%
_{\bar{t}}\in\Lambda$, we denote
\begin{align*}
\Vert\gamma_{t}\Vert:  &  =\sup\limits_{s\in\lbrack0,t]}|\gamma_{t}(s)|,\\
\Vert\gamma_{t}-\bar{\gamma}_{\bar{t}}\Vert:  &  =\sup\limits_{s\in
\lbrack0,t\vee\bar{t}]}|\gamma_{t}(s\wedge t)-\bar{\gamma}_{\bar{t}}%
(s\wedge\bar{t})|,\\
d_{\infty}(\gamma_{t},\bar{\gamma}_{\bar{t}}):  &  =\sup_{0\leq s\leq
t\vee\bar{t}}|\gamma_{t}(s\wedge t)-\bar{\gamma}_{\bar{t}}(s\wedge\bar
{t})|+|t-\bar{t}|.
\end{align*}
It is obvious that $\Lambda_{t}$ is a Banach space with respect to $\Vert
\cdot\Vert$ and $d_{\infty}$ is not a norm.

\begin{definition}
A function $u:\Lambda\mapsto\mathbb{R}$ is said to be $\Lambda$--continuous at
$\gamma_{t}\in\Lambda$, if for any $\varepsilon>0$ there exists $\delta>0$
such that for each $\bar{\gamma}_{\bar{t}}\in\Lambda$ with $d_{\infty}%
(\gamma_{t},\bar{\gamma}_{\bar{t}})<\delta$, we have $|u(\gamma_{t}%
)-u(\bar{\gamma}_{\bar{t}})|<\varepsilon$. $u$ is said to be $\Lambda
$--continuous if it is $\Lambda$--continuous at each $\gamma_{t}\in\Lambda$.
\end{definition}

\begin{definition}
Let $u:\Lambda\mapsto\mathbb{R}$ and $\gamma_{t}\in\Lambda$ be given. If there
exists $p\in\mathbb{R}^{d}$, such that
\[
u(\gamma_{t}^{x})=u(\gamma_{t})+\langle p,x\rangle+o(|x|)\ \text{as}%
\ x\rightarrow0, \ x\in\mathbb{R}^{d}.\ \
\]
Then we say that $u$ is (vertically) differentiable at $\gamma_{t}$ and denote
the gradient of $D_{x}u(\gamma_{t})=p$. $u$ is said to be vertically
differentiable in $\Lambda$ if $D_{x}u(\gamma_{t})$ exists for each
$\gamma_{t}\in\Lambda$. We can similarly define the Hessian $D_{xx}%
u(\gamma_{t})$. It is an $\mathbb{S}(d)$-valued function defined on $\Lambda$,
where $\mathbb{S}(d)$ is the space of all $d\times d$ symmetric matrices.
\end{definition}

For each $\gamma_{t}\in\Lambda$ we denote
\[
\gamma_{t,s}(r)=\gamma_{t}(r)\mathbf{1}_{[0,t)}(r)+\gamma_{t}(t)\mathbf{1}%
_{[t,s]}(r),\ \ r\in\lbrack0,s].
\]
It is clear that $\gamma_{t,s}\in\Lambda_{s}$.

\begin{definition}
For a given $\gamma_{t}\in\Lambda$ if we have
\[
u(\gamma_{t,s})=u(\gamma_{t})+a(s-t)+o(|s-t|)\ \text{as}\ s\rightarrow t,
\ s\geq t,\ \
\]
then we say that $u(\gamma_{t})$ is (horizontally) differentiable in $t$ at
$\gamma_{t}$ and denote $D_{t}u(\gamma_{t})=a$. $u$ is said to be horizontally
differentiable in $\Lambda$ if $D_{t}u(\gamma_{t})$ exists for each
$\gamma_{t}\in\Lambda$.
\end{definition}

\begin{definition}
Define $\mathbb{C}^{j,k}(\Lambda)$ as the set of function $u:=(u(\gamma
_{t}))_{\gamma_{t}\in\Lambda}$ defined on $\Lambda$ which are $j$ times
horizontally and $k$ times vertically differentiable in $\Lambda$ such that
all these derivatives are $\Lambda$--continuous.
\end{definition}

The following It\^{o} formula was firstly obtained by Dupire \cite{Dupire.B}
and then generalized by Cont and Fourni\'{e} \cite{Cont.R}, \cite{Cont-2} and
\cite{Cont-3}.

\begin{theorem}
[Functional It\^{o}'s formula]\label{w2} Let $(\Omega,\mathcal{F}%
,(\mathcal{F}_{t})_{t\in\lbrack0,T]},P)$ be a probability space, if $X$ is a
continuous semi-martingale and $u$ is in $\mathbb{C}^{1,2}(\Lambda)$, then for
any $t\in\lbrack0,T)$,
\[%
\begin{split}
u(X_{t})-u(X_{0})=  &  \int_{0}^{t}D_{s}u(X_{s})\,ds+\int_{0}^{t}D_{x}%
u(X_{s})\,dX(s)\\
&  +\frac{1}{2}\int_{0}^{t}D_{xx}u(X_{s})\,d\langle X\rangle(s)\quad
\quad\ P-a.s.
\end{split}
\]

\end{theorem}

\subsection{Non-Markovian fully coupled FBSDEs}

Let $\Omega=C([0,T];\mathbb{R}^{d})$ and $P$ the Wiener measure on
$(\Omega,\mathbb{B}(\Omega))$. We denote by $W=(W(t)_{t\in\lbrack0,T]})$ the
cannonical Wiener process, with $W(t,\omega)=\omega(t)$, $t\in\lbrack0,T]$,
$\omega\in\Omega$. For any $t\in\lbrack0,T]$. we denote by $\mathcal{F}_{t}$
the $P$-completion of $\sigma(W(s),s\in\lbrack0,t])$.

For any $t\in\lbrack0,T]$, we denote by $L^{2}(\Omega,\mathcal{F}%
_{t};\mathbb{R}^{n})$ the set of all square integrable $\mathcal{F}_{t}%
-$measurable random variables, $M^{2}(0,T;\mathbb{R}^{n})$ the set of all
$\mathbb{R}^{n}$-valued $\mathcal{F}_{t}$-adapted processes $\vartheta(\cdot)$
such that
\[
E\int_{0}^{T}\mid\vartheta(s)\mid^{2}ds<+\infty.
\]

Let $t\in\lbrack0,T]$, $\gamma_{t}\in{\Lambda}$ and $x\in\mathbb{R}^{n}$. We
consider the following Non-Markovian forward-backward SDEs:
\begin{equation}
X^{\gamma_{t},x}(s)=x+\int_{t}^{s}b(W_{r}^{\gamma_{t}},X^{\gamma_{t}%
,x}(r),Y^{\gamma_{t},x}(r),Z^{\gamma_{t},x}(r))dr+\int_{t}^{s}\sigma
(W_{r}^{\gamma_{t}},X^{\gamma_{t},x}(r),Y^{\gamma_{t},x}(r),Z^{\gamma_{t}%
,x}(r))dW(r), \tag{2.1}%
\end{equation}%
\begin{equation}
Y^{\gamma_{t},x}(s)=g(W_{T}^{\gamma_{t}},X^{\gamma_{t},x}(T))-\int_{s}%
^{T}h(W_{r}^{\gamma_{t}},X^{\gamma_{t},x}(r),Y^{\gamma_{t},x}(r),Z^{\gamma
_{t},x}(r))dr-\int_{s}^{T}Z^{\gamma_{t},x}(r)dW(r), \tag{2.2}%
\end{equation}
for every $s\in\lbrack t,T],P-a.s.$, where%
\[
W^{\gamma_{t}}(s):=\gamma(s)I_{0\leq s<t}+(W(s)-W(t)+\gamma(t))I_{t\leq s\leq
T}.
\]
We suppose that $X(t):=xI_{0\leq r\leq t}$; the processes $X,Y,Z$ take values
in $\mathbb{R}^{n},\mathbb{R}^{n},\mathbb{R}^{n\times d}$; $b,h,\sigma$ and
$g$ take values in $\mathbb{R}^{n},\mathbb{R}^{n},\mathbb{R}^{n\times d}$ and
$\mathbb{R}^{n}$. $(2.1)$ and $(2.2)$ can be rewriten as:
\begin{align*}
&  dX^{\gamma_{t},x}(s)=b(W_{s}^{\gamma_{t}},X^{\gamma_{t},x}(s),Y^{\gamma
_{t},x}(s),Z^{\gamma_{t},x}(s))ds+\sigma(W_{s}^{\gamma_{t}},X^{\gamma_{t}%
,x}(s),Y^{\gamma_{t},x}(s),Z^{\gamma_{t},x}(s))dW(s),\\
&  dY^{\gamma_{t},x}(s)=h(W_{s}^{\gamma_{t}},X^{\gamma_{t},x}(s),Y^{\gamma
_{t},x}(s),Z^{\gamma_{t},x}(s))ds+Z^{\gamma_{t},x}(s)dW(s),\\
&  X^{\gamma_{t},x}(t)=x,\qquad Y^{\gamma_{t},x}(T)=g(W_{T}^{\gamma_{t}%
},X^{\gamma_{t},x}(T)).
\end{align*}

For $z\in\mathbb{R}^{n\times d}$, we define $|z|=\{tr(zz^{T})\}^{1/2}$. For
$z^{1}\in\mathbb{R}^{n\times d}$, $z^{2}\in\mathbb{R}^{n\times d}$,
\[
((z^{1},z^{2}))=tr(z^{1}(z^{2})^{T}),
\]
and for $u^{1}=(x^{1},y^{1},z^{1})\in\mathbb{R}^{n}\times\mathbb{R}^{n}%
\times\mathbb{R}^{n\times d}$, $u^{2}=(x^{2},y^{2},z^{2})\in\mathbb{R}%
^{n}\times\mathbb{R}^{n}\times\mathbb{R}^{n\times d}$.
\[
\lbrack u^{1},u^{2}]=\langle x^{1},x^{2}\rangle+\langle y^{1},y^{2}%
\rangle+((z^{1},z^{2})).
\]
For $u=(x,y,z)\in\mathbb{R}^{n}\times\mathbb{R}^{n}\times\mathbb{R}^{n\times
d}$ and $\gamma_{t}\in{\Lambda}^{d},$ denote
\[
f(\gamma_{t},u)=(h(\gamma_{t},u),b(\gamma_{t},u),\sigma(\gamma_{t},u)).
\]

We give the following assumptions:

\begin{assumption}
For each $u\in\mathbb{R}^{n}\times\mathbb{R}^{n}\times\mathbb{R}^{n\times d}$,
$f(\cdot,u)\in M^{2}(0,T;\mathbb{R}^{n}\times\mathbb{R}^{n}\times
\mathbb{R}^{n\times d})$, and for each $x\in\mathbb{R}^{n}$, $g(\cdot,x)\in
L^{2}(\Omega,\mathcal{F}_{T};\mathbb{R}^{n})$; there exists a constant
$c_{1}>0$, such that%
\[%
\begin{array}
[c]{c}%
\mid f(\gamma_{t},u^{1})-f(\gamma_{t},u^{2})\mid\leq c_{1}\mid u^{1}-u^{2}%
\mid,\;a.e.t\in\lbrack0,T],\\
\forall\gamma_{t}\in{\Lambda}\text{ and }u^{1}\in\mathbb{R}^{n}\times
\mathbb{R}^{n}\times\mathbb{R}^{n\times d},u^{2}\in\mathbb{R}^{n}%
\times\mathbb{R}^{n}\times\mathbb{R}^{n\times d};
\end{array}
\]
and
\[
\mid g(\gamma_{t},x^{1})-g(\gamma_{t},x^{2})\mid\leq c_{1}\mid x^{1}-x^{2}%
\mid,\quad\forall\gamma_{t}\in{\Lambda}\text{ and }(x^{1},x^{2})\in
\mathbb{R}^{n}\times\mathbb{R}^{n}.
\]

\end{assumption}

\begin{assumption}
There exists a constant $c_{2}$%
$>$%
0, such that%
\[%
\begin{array}
[c]{c}%
\lbrack f(\gamma_{t},u^{1})-f(\gamma_{t},u^{2}),u^{1}-u^{2})]\leq-c_{2}\mid
u^{1}-u^{2}\mid^{2},\quad a.e.t\in\lbrack0,T],\\
\forall\gamma_{t}\in{\Lambda}\text{ and }u^{1}\in\mathbb{R}^{n}\times
\mathbb{R}^{n}\times\mathbb{R}^{n\times d},u^{2}\in\mathbb{R}^{n}%
\times\mathbb{R}^{n}\times\mathbb{R}^{n\times d};
\end{array}
\]
\newline and
\[
\langle g(\gamma_{t},x^{1})-g(\gamma_{t},x^{2}),x^{1}-x^{2}\rangle\geq
c_{2}^{2}\mid x^{1}-x^{2}\mid,\quad\forall\gamma_{t}\in{\Lambda}\text{ and
}(x^{1},x^{2})\in\mathbb{R}^{n}\times\mathbb{R}^{n}.
\]

\end{assumption}

\begin{definition}
A triple $(X,Y,Z)$:$[0,T]\times\Omega$$\rightarrow$ $\mathbb{R}^{n}%
\times\mathbb{R}^{n}\times\mathbb{R}^{n\times d}$ is called an adapted
solution of the Eqs.$(2.1)$ and $(2.2)$, if $(X,Y,Z)\in M^{2}(0,T;\mathbb{R}%
^{n}\times\mathbb{R}^{n}\times\mathbb{R}^{n\times d})$, and it satisfies
$(2.1)$ and $(2.2)$ $P-a.s.$
\end{definition}

Then we have the following theorem (Theorem 3.1 in \cite{Hu Y}):

\begin{theorem}
Let Assumptions $2.1$ and $2.2$ hold, then there exists a unique adapted
solution $(X,Y,Z)$ for Eqs. $(2.1)$ and $(2.2)$.
\end{theorem}

For the theory of FBSDEs, the readers may refer to \cite{Pardoux.E}, \cite{Ma
J}, \cite{Peng S} and the references in \cite{Ma J-1}.

\section{Regularity}

We first recall some notions in Pardoux and Peng \cite{Pardoux.E 2}.
$C^{n}(\mathbb{R}^{p};\mathbb{R}^{q})$, $C_{b}^{n}(\mathbb{R}^{p}%
;\mathbb{R}^{q})$, $C_{p}^{n}(\mathbb{R}^{p};\mathbb{R}^{q})$ will denote
respectively the set of functions of class $C^{n}$ from $\mathbb{R}^{p}$ into
$\mathbb{R}^{q}$, the set of those functions of class $C_{b}^{n}$ whose
partial derivatives of order less than or equal to $n$ are bounded, and the
set of those functions of class $C_{p}^{n}$ which, together with all their
partial derivatives of order less than or equal to n, grow at most like a
polynomial function of the variable $x$ at infinity.

We give the definition of derivatives in our context.

\begin{definition}
An $\mathbb{R}^{n}$-valued function $g(\gamma,x)$ on $\Lambda_{T}%
\times\mathbb{R}^{n}$ is said to be in $C^{2,2}(\Lambda_{T}\times
\mathbb{R}^{n})$, if there exists the second order partial derivative of $g$
in $x$ and for each ${\gamma_{t}},\gamma_{\gamma_{t}^{y}}\in\Lambda_{T}$,
$t\in\lbrack0,T]$, there exist $p_{1}\in\mathbb{R}^{d}$ and $p_{2}%
\in\mathbb{R}^{d}\times\mathbb{R}^{d}$ such that $p_{2}$ is symmetric and
\[
g(\gamma_{\gamma_{t}^{y}},x)-g({\gamma_{t}},x)=\langle p_{1},y\rangle+\frac
{1}{2}\langle p_{2}y,y\rangle+o(|y|^{2}),\quad y\in\mathbb{R}^{d},
\]
where $\gamma_{\gamma_{t}^{y}}=\gamma(r)I_{[0,t)}(r)+(\gamma(r)+y)I_{[t,T]}%
(r)$. We denote
\[
g{^{\prime}}_{\gamma_{t}}({\gamma_{t},x}):=p_{1},\text{and }g{^{\prime\prime}%
}_{\gamma_{t}}({\gamma_{t},x}):=p_{2}.
\]
$g$ is said to be in $C_{l,lip}^{2}(\Lambda_{T})$ if $g{^{\prime}}_{\gamma
_{t}}(\gamma)$ and $g{^{\prime\prime}}_{\gamma_{t}}(\gamma)$ exist for each
$\gamma\in\Lambda_{T}$, and there exists some constants $C\geq0$ and $k\geq0$
depending only on $g$ such that for each $\gamma,\bar{\gamma}\in\Lambda
_{T},t,s\in\lbrack0,T]$,
\begin{align*}
&  \mid g(\gamma)-g(\bar{\gamma})\mid\leq C(\parallel\gamma\parallel
^{k}+\parallel\bar{\gamma}\parallel^{k})\parallel\gamma-\bar{\gamma}%
\parallel,\\
&  \mid\Phi_{\gamma_{t}}(\gamma)-\Phi_{\gamma_{s}}(\bar{\gamma})\mid\leq
C(\parallel\gamma\parallel^{k}+\parallel\bar{\gamma}\parallel^{k})(\mid
t-s\mid+\parallel\gamma-\bar{\gamma}\parallel)
\end{align*}
with $\Phi=g{^{\prime}}_{\gamma_{t}}(\gamma),g{^{\prime\prime}}_{\gamma_{t}%
}(\gamma)$. Analogously, we can define $C^{2}(\Lambda_{t})$, $C_{l,lip}%
^{2}(\Lambda_{t})$, $C_{l,lip}^{1}(\Lambda_{t})$, $C_{l,lip}(\Lambda_{t})$.
\end{definition}

Now we reconsider the solvability of equation (2.1) and (2.2). We use $\Psi$
to represent $b,\sigma$or $h$.

\begin{assumption}
$g\in C_{l,lip}^{2,2}(\Lambda_{T}\times\mathbb{R}^{n})$ with the Lipschitz
constants $C$,$k$ and the first order partial derivatives in x are bounded, as
well as their derivatives of order one and two with respect to x..
\end{assumption}

\begin{assumption}
Let $\Psi(\gamma_{t},x,y,z)=\bar{\Psi}(t,\gamma(t),x,y,z)$, where $\bar{\Psi
}:[0,T]\times\mathbb{R}^{d}\times\mathbb{R}^{n}\times\mathbb{R}^{n}%
\times\mathbb{R}^{n\times d}\mapsto\mathbb{R}^{n}$. Suppose that
$(t,r,x,y,z)\mapsto\bar{\Psi}(t,r,x,y,z)$ is of class $C_{p}^{0,3}%
([0,T]\times\mathbb{R}^{d}\times\mathbb{R}^{n}\times\mathbb{R}^{n}%
\times{\mathbb{R}^{n\times d}};\mathbb{R}^{n})$ and the first order partial
derivatives in r,x,y and z are bounded, as well as their derivatives of up to
order two with repect to x,y,z.
\end{assumption}

By Theorem 2.1, it is easy to see that under Assumptions 2.2, 3.1 and 3.2, the
FBSDE (2.1) and (2.2) has a uniqueness solution.

\subsection{Regularity of the solution of FBSDEs}

We first fix $x$ and assume that the Lipschitz constants with respect to
$\Psi$ are C and k. Then we establish second order moment estimates for the
solution of FBSDE (2.1) and (2.2).

\begin{lemma}
Under Assumptions 2.2, 3.1 and 3.2, there exists $C_{2}$ and $q$ depending
only on $C,T,k,x$ such that
\begin{align*}
&  E[\sup_{s\in\lbrack t,T]}\mid X^{\gamma_{t},x}(s)\mid^{2}]\leq
C_{2}(1+\parallel\gamma_{t}\parallel^{q}),\\
&  E[\sup_{s\in\lbrack t,T]}\mid Y^{\gamma_{t},x}(s)\mid^{2}]\leq
C_{2}(1+\parallel\gamma_{t}\parallel^{q}),\\
&  E[(\int_{t}^{T}\mid Z^{\gamma_{t},x}(s)\mid^{2}ds)]\leq C_{2}%
(1+\parallel\gamma_{t}\parallel^{q}).
\end{align*}

\end{lemma}

\textbf{\noindent Proof. }For simplicity, we only study the case
$n=d=1$.\newline Applying It\^{o}'s formula to $(Y_{\gamma_{t},x}%
(s))^{2}e^{\beta_{1}s}$ yields that%
\[%
\begin{array}
[c]{rl}
& (Y^{\gamma_{t},x}(s))^{2}e^{\beta_{1}s}+\int_{s}^{T}e^{\beta_{1}%
r}[(Z^{\gamma_{t},x}(r))^{2}+\beta_{1}(Y^{\gamma_{t},x}(r))^{2}]dr\\
= & g^{2}(W_{T}^{\gamma_{t}},X^{\gamma_{t},x}(T))e^{\beta_{1}T}-\int_{s}%
^{T}2e^{\beta_{1}r}Y^{\gamma_{t},x}(r)h(W_{r}^{\gamma_{t}},X^{\gamma_{t}%
,x}(r),Y^{\gamma_{t},x}(r),Z^{\gamma_{t},x}(r))dr\\
& -\int_{s}^{T}2e^{\beta_{1}r}Y^{\gamma_{t},x}(r)Z^{\gamma_{t},x}(r)dW(r).
\end{array}
\]
So%
\[%
\begin{array}
[c]{rl}
& (Y^{\gamma_{t},x}(s))^{2}+E[\int_{s}^{T}e^{\beta_{1}(r-s)}[(Z^{\gamma_{t}%
,x}(r))^{2}+\beta_{1}(Y^{\gamma_{t},x}(r))^{2}]dr\mid\mathcal{F}_{s}]\\
= & E[g^{2}(W_{T}^{\gamma_{t}},X^{\gamma_{t},x}(T))e^{\beta_{1}(T-s)}%
\mid\mathcal{F}_{s}]\\
& -E[\int_{s}^{T}2e^{\beta_{1}(r-s)}Y^{\gamma_{t},x}(r)h(W_{r}^{\gamma_{t}%
},X^{\gamma_{t},x}(r),Y^{\gamma_{t},x}(r),Z^{\gamma_{t},x}(r))dr\mid
\mathcal{F}_{s}].
\end{array}
\]
Then we have%
\[%
\begin{array}
[c]{rl}
& E\sup_{t\leq s\leq T}(Y^{\gamma_{t},x}(s))^{2}+E[\int_{t}^{T}e^{\beta
_{1}(r-t)}[(Z^{\gamma_{t},x}(r))^{2}+\beta_{1}(Y^{\gamma_{t},x}(T))^{2}]dr]\\
\leq & E[g^{2}(W_{T}^{\gamma_{t}},X^{\gamma_{t},x}(r))e^{\beta_{1}%
(T-t)}]+E[\int_{t}^{T}e^{\beta_{1}(r-t)}\frac{2}{\beta_{1}}h^{2}(W_{r}%
^{\gamma_{t}},X^{\gamma_{t},x}(r),Y^{\gamma_{t},x}(r),Z^{\gamma_{t}%
,x}(r))dr]\\
& +E[\int_{t}^{T}e^{\beta_{1}(r-t)}\frac{\beta_{1}}{2}(Y^{\gamma_{t}%
,x}(r))^{2}d(r)]
\end{array}
\]
and%
\begin{equation}%
\begin{array}
[c]{rl}
& E\sup_{t\leq s\leq T}(Y^{\gamma_{t},x}(s))^{2}+E[\int_{t}^{T}e^{\beta
_{1}(r-t)}[(Z^{\gamma_{t},x}(r))^{2}+\frac{\beta_{1}}{2}(Y^{\gamma_{t}%
,x}(r))^{2}]dr]\\
\leq & E[g^{2}(W_{T}^{\gamma_{t}},X^{\gamma_{t},x}(T))e^{\beta_{1}%
(T-t)}]+E[\int_{t}^{T}e^{\beta_{1}(r-t)}\frac{2}{\beta_{1}}h^{2}(W_{r}%
^{\gamma_{t}},X^{\gamma_{t},x}(r),Y^{\gamma_{t},x}(r),Z^{\gamma_{t},x}(r))dr].
\end{array}
\tag{3.1}%
\end{equation}
Applying It\^{o}'s formula to $(X^{\gamma_{t},x}(s))^{2}$ yields that%
\[%
\begin{array}
[c]{rl}
& (X^{\gamma_{t},x}(s))^{2}\\
= & x^{2}+\int_{t}^{s}2X^{\gamma_{t},x}(r)b(W_{r}^{\gamma_{t}},X^{\gamma
_{t},x}(r),Y^{\gamma_{t},x}(r),Z^{\gamma_{t},x}(r))dr+\int_{t}^{s}\sigma
^{2}(W_{r}^{\gamma_{t}},X^{\gamma_{t},x}(r),Y^{\gamma_{t},x}(r),Z^{\gamma
_{t},x}(r))dr\\
& +\int_{t}^{s}2X^{\gamma_{t},x}(r)\sigma(W_{r}^{\gamma_{t}},X^{\gamma_{t}%
,x}(r),Y^{\gamma_{t},x}(r),Z^{\gamma_{t},x}(r))dW(r).
\end{array}
\]
By inequality $2ab\leq a^{2}+b^{2}$ and Burkholder-Davis-Gundy's inequality,
there is a constant $C_{0}$ such that%
\begin{equation}%
\begin{array}
[c]{rl}
& E\sup_{t\leq s\leq T}(X^{\gamma_{t},x}(s))^{2}\\
\leq & C_{0}[x^{2}+E\int_{t}^{T}(X^{\gamma_{t},x}(s))^{2}ds+E\int_{t}^{T}%
b^{2}(W_{s}^{\gamma_{t}},X^{\gamma_{t},x}(s),Y^{\gamma_{t},x}(s),Z^{\gamma
_{t},x}(s))ds\\
& +E\int_{t}^{T}\sigma^{2}(W_{s}^{\gamma_{t}},X^{\gamma_{t},x}(s),Y^{\gamma
_{t},x}(s),Z^{\gamma_{t},x}(s))ds].
\end{array}
\tag{3.2}%
\end{equation}
Applying It\^{o} formula to $X^{\gamma_{t},x}(s)Y^{\gamma_{t},x}(s)$,%
\[%
\begin{array}
[c]{rl}
& X^{\gamma_{t},x}(T)Y^{\gamma_{t},x}(T)-X^{\gamma_{t},x}(t)Y^{\gamma_{t}%
,x}(t)\\
= & \int_{t}^{T}X^{\gamma_{t},x}(r)h(W_{r}^{\gamma_{t}},X^{\gamma_{t}%
,x}(r),Y^{\gamma_{t},x}(r),Z^{\gamma_{t},x}(r))dr\\
& +\int_{t}^{T}X^{\gamma_{t},x}(r)Z^{\gamma_{t},x}(r)dW(r)\\
& +\int_{t}^{T}Y^{\gamma_{t},x}(r)b(W_{r}^{\gamma_{t}},X^{\gamma_{t}%
,x}(r),Y^{\gamma_{t},x}(r),Z^{\gamma_{t},x}(r))dr\\
& +\int_{t}^{T}Z^{\gamma_{t},x}(r)\sigma(W_{r}^{\gamma_{t}},X^{\gamma_{t}%
,x}(r),Y^{\gamma_{t},x}(r),Z^{\gamma_{t},x}(r))dr\\
& +\int_{t}^{T}Y^{\gamma_{t},x}(r)\sigma(W_{r}^{\gamma_{t}},X^{\gamma_{t}%
,x}(r),Y^{\gamma_{t},x}(r),Z^{\gamma_{t},x}(r))dW(r).
\end{array}
\]
Set
\[
0^{\gamma_{t},x}(s)=\gamma_{t}(s)I_{0\leq s\leq t}+0I_{t<s\leq T}(s).
\]
We have%
\[%
\begin{array}
[c]{rl}
& (X^{\gamma_{t},x}(T)-0^{\gamma_{t},x}(T))(g(W_{T}^{\gamma_{t}},X^{\gamma
_{t},x}(T))-g(W_{T}^{\gamma_{t}},0^{\gamma_{t},x}(T)))\\
& +(X^{\gamma_{t},x}(T)-0^{\gamma_{t},x}(T))g(W_{T}^{\gamma_{t}}%
,0)-X^{\gamma_{t},x}(t)Y^{\gamma_{t},x}(t)\\
= & \int_{t}^{T}[(X^{\gamma_{t},x}(r)-0^{\gamma_{t},x}(r))(h(W_{r}^{\gamma
_{t}},X^{\gamma_{t},x}(r),Y^{\gamma_{t},x}(r),Z^{\gamma_{t},x}(r))-h(W_{r}%
^{\gamma_{t}},0^{\gamma_{t},x}(r),0^{\gamma_{t},x}(r),0^{\gamma_{t},x}(r)))\\
& +(X^{\gamma_{t},x}(r)-0^{\gamma_{t},x}(r))h(W_{r}^{\gamma_{t}},0^{\gamma
_{t},x}(r),0^{\gamma_{t},x}(r),0^{\gamma_{t},x}(r))]dr\\
& +\int_{t}^{T}Y^{\gamma_{t},x}(r)\sigma(W_{r}^{\gamma_{t}},X^{\gamma_{t}%
,x}(r),Y^{\gamma_{t},x}(r),Z^{\gamma_{t},x}(r))dW(r)+\int_{t}^{T}X^{\gamma
_{t},x}(r)Z^{\gamma_{t},x}(r)dW(r)\\
& +\int_{t}^{T}[(Y^{\gamma_{t},x}(r)-0^{\gamma_{t},x}(r))(b(W_{r}^{\gamma_{t}%
},X^{\gamma_{t},x}(r),Y^{\gamma_{t},x}(r),Z^{\gamma_{t},x}(r))-b(W_{r}%
^{\gamma_{t}},0^{\gamma_{t},x}(r),0^{\gamma_{t},x}(r),0^{\gamma_{t},x}(r)))\\
& +(Y^{\gamma_{t},x}(r)-0^{\gamma_{t},x}(r))b(W_{r}^{\gamma_{t}},0^{\gamma
_{t},x}(r),0^{\gamma_{t},x}(r),0^{\gamma_{t},x}(r))]dr\\
& +\int_{t}^{T}[(Z^{\gamma_{t},x}(r)-0^{\gamma_{t},x}(r))(\sigma(W_{r}%
^{\gamma_{t}},X^{\gamma_{t},x}(r),Y^{\gamma_{t},x}(r),Z^{\gamma_{t}%
,x}(r))-\sigma(W_{r}^{\gamma_{t}},0^{\gamma_{t},x}(r),0^{\gamma_{t}%
,x}(r),0^{\gamma_{t},x}(r)))\\
& +(Z^{\gamma_{t},x}(r)-0^{\gamma_{t},x}(r))\sigma(W_{r}^{\gamma_{t}%
},0^{\gamma_{t},x}(r),0^{\gamma_{t},x}(r),0^{\gamma_{t},x}(r))]dr.
\end{array}
\]
Then by the Assumption $2.1$, for enough small $\varepsilon>0,$%
\begin{equation}%
\begin{array}
[c]{rl}
& (C_{2}-{\varepsilon})E\{(X^{\gamma_{t},x}(T))^{2}+\int_{t}^{T}%
[(X^{\gamma_{t},x}(r))^{2}+(Y^{\gamma_{t},x}(r))^{2}+(Z^{\gamma_{t},x}%
(r))^{2}]\}\\
\leq & C_{\varepsilon}+\frac{x}{\beta_{3}}+{\beta_{3}}{Y^{\gamma_{t},x}(t)}.
\end{array}
\tag{3.3}%
\end{equation}
Taking $\beta_{1}=4C^{2}+1$, and $q=2(1+k)$, from $(3.1)$, $(3.3)$ and
$(3.2)$, $(3.3)$, we derive (note that $C_{\varepsilon}$ will change line by
line) \newline%
\begin{equation}%
\begin{array}
[c]{rl}
& E\sup_{t\leq s\leq T}(Y^{\gamma_{t},x}(s))^{2}+E[\int_{t}^{T}[(Z^{\gamma
_{t},x}(r))^{2}+(Y^{\gamma_{t},x}(r))^{2}]dr]\\
\leq & C_{\varepsilon}(1+\parallel\gamma_{t}\parallel^{q}+\frac{x}{\beta_{3}%
}+{\beta_{3}}{Y^{\gamma_{t},x}(t)}).
\end{array}
\tag{3.4}%
\end{equation}
We also have
\begin{equation}
E\sup_{t\leq s\leq T}(X^{\gamma_{t},x}(s))^{2}\leq C_{\varepsilon}(1+\frac
{x}{\beta_{3}}+\parallel\gamma_{t}\parallel^{q}+{\beta_{3}}{Y^{\gamma_{t}%
,x}(t)}). \tag{3.5}%
\end{equation}
By $(3.4)$ and $(3.5)$, we have
\[%
\begin{array}
[c]{rl}
& E[\sup_{t\leq s\leq T}(Y^{\gamma_{t},x}(s))^{2}+\sup_{t\leq s\leq
T}(X^{\gamma_{t},x}(s))^{2}]+E[\int_{t}^{T}[(Z^{\gamma_{t},x}(r))^{2}%
+(Y^{\gamma_{t},x}(r))^{2}]dr]\\
\leq & C_{\varepsilon}(1+\parallel\gamma_{t}\parallel^{q}+\frac{x}{\beta_{3}%
}+{\beta_{3}}{Y^{\gamma_{t},x}(t)}).
\end{array}
\]
Finally taking $\beta_{3}$ enough small, we get the result. \quad$\Box$

Now we study the regularity properties of the solution of FBSDE (2.1) and
(2.2) with respect to the "parameter" $\gamma_{t}$. For $0\leq s<t\leq T$, set
$Y^{\gamma_{t},x}(s)=Y^{\gamma_{t},x}(s\vee t)$ and $Z^{\gamma_{t},x}(s)=0$.

\begin{theorem}
Under Assumptions $2.2$, $3.2$, $3.3$,there exists $C_{2}$ and $q$ depending
only on $C,c_{2},x$ such that for any $t,\bar{t}\in\lbrack0,T]$, $\gamma
_{t},\bar{\gamma}_{\bar{t}}$, and $h,\bar{h}\in{\mathbb{R}\setminus\{0\}}$.
\end{theorem}

\noindent$(i)$
\[
E[\sup_{u\in\lbrack t\vee\bar{t},T]}\mid Y^{\gamma_{t},x}(u)-Y^{\bar{\gamma
}_{\bar{t}},x}(u)\mid^{2}]\leq C_{2}(1+\parallel\gamma_{t}\parallel
^{q}+\parallel\bar{\gamma}_{\bar{t}}\parallel^{q})(\parallel\gamma_{t}%
-\bar{\gamma}_{\bar{t}}\parallel^{2}+\mid t-\bar{t}\mid),
\]
$(ii)$
\[
E[\sup_{u\in\lbrack t\vee\bar{t},T]}\mid X^{\gamma_{t},x}(u)-X^{\bar{\gamma
}_{\bar{t}},x}(u)\mid^{2}]\leq C_{2}(1+\parallel\gamma_{t}\parallel
^{q}+\parallel\bar{\gamma}_{\bar{t}}\parallel^{q})(\parallel\gamma_{t}%
-\bar{\gamma}_{\bar{t}}\parallel^{2}+\mid t-\bar{t}\mid),
\]
$(iii)$
\[
E[\int_{t\vee\bar{t}}^{T}\mid Z^{\gamma_{t},x}(u)-Z^{\bar{\gamma}_{\bar{t}}%
,x}(u)\mid^{2}du]\leq C_{2}(1+\parallel\gamma_{t}\parallel^{q}+\parallel
\bar{\gamma}_{\bar{t}}\parallel^{q})(\parallel\gamma_{t}-\bar{\gamma}_{\bar
{t}}\parallel^{2}+\mid t-\bar{t}\mid),
\]
$(iv)$
\[%
\begin{array}
[c]{rl}
& E[\sup_{u\in\lbrack t\vee\bar{t},T]}\mid\Delta_{h}^{i}Y^{\gamma_{t}%
,x}(u)-\Delta_{h}^{i}Y^{\bar{\gamma}_{{\bar{t}}},x}(u)\mid^{2}]\\
\leq & C_{2}(1+\parallel\gamma_{t}\parallel^{q}+\parallel\bar{\gamma}_{\bar
{t}}\parallel^{q}+\mid h\mid^{q}+\mid\bar{h}\mid^{q}))(\mid h-\bar{h}\mid
^{2}+\parallel\gamma_{t}-\bar{\gamma}_{\bar{t}}\parallel^{2}+\mid t-\bar
{t}\mid),
\end{array}
\]
$(v)$
\[%
\begin{array}
[c]{rl}
& E[\sup_{u\in\lbrack t\vee\bar{t},T]}\mid\Delta_{h}^{i}X^{\gamma_{t}%
,x}(u)-\Delta_{h}^{i}X^{\bar{\gamma}_{{\bar{t}}},x}(u)\mid^{2}]\\
\leq & C_{2}(1+\parallel\gamma_{t}\parallel^{q}+\parallel\bar{\gamma}_{\bar
{t}}\parallel^{q}+\mid h\mid^{q}+\mid\bar{h}\mid^{q}))(\mid h-\bar{h}\mid
^{2}+\parallel\gamma_{t}-\bar{\gamma}_{\bar{t}}\parallel^{2}+\mid t-\bar
{t}\mid),
\end{array}
\]
$(vi)$
\[%
\begin{array}
[c]{rl}
& E[\int_{t\vee\bar{t}}^{T}\mid\Delta_{h}^{i}Z^{\gamma_{t},x}(u)-\Delta
_{h}^{i}Z^{\bar{\gamma}_{{\bar{t}}},x}(u)\mid^{2}du]\\
\leq & C_{2}(1+\parallel\gamma_{t}\parallel^{q}+\parallel\bar{\gamma}_{\bar
{t}}\parallel^{q}+\mid h\mid^{q}+\mid\bar{h}\mid^{q}))(\mid h-\bar{h}\mid
^{2}+\parallel\gamma_{t}-\bar{\gamma}_{\bar{t}}\parallel^{2}+\mid t-\bar
{t}\mid),
\end{array}
\]
where%
\[%
\begin{array}
[c]{rl}%
\Delta_{h}^{i}X^{\gamma_{t},x}(s)= & \frac{1}{h}(X^{\gamma_{t}^{h_{e_{i}}}%
,x}(s)-X^{\gamma_{t},x}(s)),\\
\Delta_{h}^{i}Y^{\gamma_{t},x}(s)= & \frac{1}{h}(Y^{\gamma_{t}^{h_{e_{i}}}%
,x}(s)-Y^{\gamma_{t},x}(s)),\Delta_{h}^{i}Z^{\gamma_{t},x}(s)=\frac{1}%
{h}(Z^{\gamma_{t}^{h_{e_{i}}},x}(s)-Z^{\gamma_{t},x}(s))
\end{array}
\]
and $(e_{1},\cdots,e_{d})$ is an orthonormal basis of $\mathbb{R}^{d}$.

$\noindent$\textbf{{Proof.}} $(X^{\gamma_{t},x}-X^{\bar{\gamma}_{\bar{t}}%
,x},Y^{\gamma_{t},x}-Y^{\bar{\gamma}_{\bar{t}},x},Z^{\gamma_{t},x}%
-Z^{\bar{\gamma}_{\bar{t}},x})$ can be formed as a linearized FBSDE: for each
$s\in\lbrack{t},T]$ and $\bar{t}\leq{t}$,%
\[%
\begin{array}
[c]{rl}
& X^{\gamma_{t},x}(s)-X^{\bar{\gamma}_{\bar{t}},x}(s)\\
= & \int_{\bar{t}}^{t}b(W_{r}^{\bar{\gamma}_{\bar{t}}},X^{\gamma_{\bar{t}}%
,x}(r),Y^{\bar{\gamma}_{\bar{t}},x}(r),Z^{\bar{\gamma}_{\bar{t}},x}(r))dr\\
& +\int_{t}^{s}b(W_{r}^{\gamma_{t}},X^{\gamma_{t},x}(r),Y^{\gamma_{t}%
,x}(r),Z^{\gamma_{t},x}(r))-b(W_{r}^{\bar{\gamma}_{\bar{t}}},X^{\gamma
_{\bar{t}},x}(r),Y^{\bar{\gamma}_{\bar{t}},x}(r),Z^{\bar{\gamma}_{\bar{t}}%
,x}(r))dr\\
& +\int_{t}^{s}\sigma(W_{r}^{\gamma_{t}},Y^{\gamma_{t},x}(r),X^{\gamma_{t}%
,x}(r),Z^{\gamma_{t},x}(r))-\sigma(W_{r}^{\bar{\gamma}_{\bar{t}}}%
,X^{\gamma_{\bar{t}},x}(r),Y^{\bar{\gamma}_{\bar{t}},x}(r),Z^{\bar{\gamma
}_{\bar{t}},x}(r))dW(r)\\
= & \int_{\bar{t}}^{t}b(W_{r}^{\bar{\gamma}_{\bar{t}}},X^{\gamma_{\bar{t}}%
,x}(r),Y^{\bar{\gamma}_{\bar{t}},x}(r),Z^{\bar{\gamma}_{\bar{t}},x}%
(r))dr+\int_{t}^{s}[\alpha_{\gamma_{t},\bar{\gamma}_{\bar{t}}}(r)+\theta
_{\gamma_{t},\bar{\gamma}_{\bar{t}}}(r)(X^{\gamma_{t},x}(r)-X^{\bar{\gamma
}_{\bar{t}},x}(r))\\
& +\beta_{\gamma_{t},\bar{\gamma}_{\bar{t}}}(r)(Y^{\gamma_{t},x}%
(r)-Y^{\bar{\gamma}_{\bar{t}},x}(r))+\delta_{\gamma_{t},\bar{\gamma}_{\bar{t}%
}}(Z^{\gamma_{t},x}(r)-Z^{\bar{\gamma}_{\bar{t}},x}(r))]dr\\
& +\int_{t}^{s}[\bar{\alpha}_{\gamma_{t},\bar{\gamma}_{\bar{t}}}%
(r)+\bar{\theta}_{\gamma_{t},\bar{\gamma}_{\bar{t}}}(X^{\gamma_{t}%
,x}(r)-X^{\bar{\gamma}_{\bar{t}},x}(r))+\bar{\beta}_{\gamma_{t},\bar{\gamma
}_{\bar{t}}}(Y^{\gamma_{t},x}(r)\\
& -Y^{\bar{\gamma}_{\bar{t}},x}(r))+\bar{\delta}_{\gamma_{t},\bar{\gamma
}_{\bar{t}}}(Z^{\gamma_{t},x}(r)-Z^{\bar{\gamma}_{\bar{t}},x}(r))]dW(r),
\end{array}
\]
and%
\[%
\begin{array}
[c]{rl}
& Y^{\gamma_{t},x}(s)-Y^{\bar{\gamma}_{\bar{t}},x}(s)\\
= & g(W_{T}^{\gamma_{t}},X^{\gamma_{t},x}(T))-g(W_{T}^{\bar{\gamma}_{\bar{t}}%
},X^{\bar{\gamma}_{\bar{t}},x}(T))+\int_{s}^{T}[h(W_{r}^{\gamma_{t}}%
,X^{\gamma_{t},x}(r),Y^{\gamma_{t},x}(r),Z^{\gamma_{t},x}(r))\\
& -h(W_{r}^{\bar{\gamma}_{\bar{t}}},X^{\gamma_{\bar{t}},x}(r),Y^{\bar{\gamma
}_{\bar{t}},x}(r),Z^{\bar{\gamma}_{\bar{t}},x}(r))]dr+\int_{s}^{T}%
(Z^{\gamma_{t},x}(r)-Z^{\bar{\gamma}_{\bar{t}},x}(r))dW(r)\\
= & g(W_{T}^{\gamma_{t}},X^{\gamma_{t},x}(T))-g(W_{T}^{\bar{\gamma}_{\bar{t}}%
},X^{\bar{\gamma}_{\bar{t}},x}(T))-\int_{s}^{T}[\hat{\alpha}_{\gamma_{t}%
,\bar{\gamma}_{\bar{t}}}(r)+\hat{\theta}_{\gamma_{t},\bar{\gamma}_{\bar{t}}%
}(X^{\gamma_{t},x}(r)-X^{\bar{\gamma}_{\bar{t}},x}(r))\\
& +\hat{\beta}_{\gamma_{t},\bar{\gamma}_{\bar{t}}}(Y^{\gamma_{t},x}%
(r)-Y^{\bar{\gamma}_{\bar{t}},x}(r))+\hat{\delta}_{\gamma_{t},\bar{\gamma
}_{\bar{t}}}(Z^{\gamma_{t},x}(r)-Z^{\bar{\gamma}_{\bar{t}},x}(r))]dr+\int%
_{s}^{T}(Z^{\gamma_{t},x}(r)-Z^{\bar{\gamma}_{\bar{t}},x}(r))]dW(r).
\end{array}
\]
Here%
\[%
\begin{array}
[c]{rl}%
\alpha_{\gamma_{t},\bar{\gamma}_{\bar{t}}}(r) & =b(W_{r}^{{\gamma}_{{t}}%
},X^{\bar{\gamma}_{\bar{t}},x}(r),Y^{\bar{\gamma}_{\bar{t}},x}(r),Z^{\bar
{\gamma}_{\bar{t}},x}(r))-b(W_{r}^{\bar{\gamma}_{\bar{t}}},X^{\bar{\gamma
}_{\bar{t}},x}(r),Y^{\bar{\gamma}_{\bar{t}},x}(r),Z^{\bar{\gamma}_{\bar{t}}%
,x}(r)),\\
\theta_{\gamma_{t},\bar{\gamma}_{\bar{t}}}(r) & =\int_{0}^{1}\frac{\partial
b}{\partial x}(W_{r}^{\gamma_{t}},U^{\bar{\gamma}_{\bar{t}},x}(r)+\theta
(U^{{\gamma}_{{t}}}(r)-U^{\bar{\gamma}_{\bar{t}},x}(r)))d\theta,\\
\beta_{\gamma_{t},\bar{\gamma}_{\bar{t}}}(r) & =\int_{0}^{1}\frac{\partial
b}{\partial y}(W_{r}^{\gamma_{t}},U^{\bar{\gamma}_{\bar{t}},x}(r)+\theta
(U^{{\gamma}_{{t}}}(r)-U^{\bar{\gamma}_{\bar{t}},x}(r)))d\theta,\\
\delta_{\gamma_{t},\bar{\gamma}_{\bar{t}}}(r) & =\int_{0}^{1}\frac{\partial
b}{\partial z}(W_{r}^{\gamma_{t}},U^{\bar{\gamma}_{\bar{t}},x}(r)+\theta
(U^{{\gamma}_{{t}}}(r)-U^{\bar{\gamma}_{\bar{t}},x}(r)))d\theta,\\
\bar{\alpha}_{\gamma_{t},\bar{\gamma}_{\bar{t}}}(r) & =\sigma(W_{r}^{{\gamma
}_{{t}}},X^{\bar{\gamma}_{\bar{t}},x}(r),Y^{\bar{\gamma}_{\bar{t}}%
,x}(r),Z^{\bar{\gamma}_{\bar{t}},x}(r))-\sigma(W_{r}^{\bar{\gamma}_{\bar{t}}%
},X^{\bar{\gamma}_{\bar{t}},x}(r),Y^{\bar{\gamma}_{\bar{t}},x}(r),Z^{\bar
{\gamma}_{\bar{t}},x}(r)),\\
\bar{\theta}_{\gamma_{t},\bar{\gamma}_{\bar{t}}}(r) & =\int_{0}^{1}%
\frac{\partial\sigma}{\partial x}(W_{r}^{\gamma_{t}},U^{\bar{\gamma}_{\bar{t}%
},x}(r)+\theta(U^{{\gamma}_{{t}}}(r)-U^{\bar{\gamma}_{\bar{t}},x}%
(r)))d\theta,\\
\bar{\beta}_{\gamma_{t},\bar{\gamma}_{\bar{t}}}(r) & =\int_{0}^{1}%
\frac{\partial\sigma}{\partial y}(W_{r}^{\gamma_{t}},U^{\bar{\gamma}_{\bar{t}%
},x}(r)+\theta(U^{{\gamma}_{{t}}}(r)-U^{\bar{\gamma}_{\bar{t}},x}%
(r)))d\theta,\\
\bar{\delta}_{\gamma_{t},\bar{\gamma}_{\bar{t}}}(r) & =\int_{0}^{1}%
\frac{\partial\sigma}{\partial z}(W_{r}^{\gamma_{t}},U^{\bar{\gamma}_{\bar{t}%
},x}(r)+\theta(U^{{\gamma}_{{t}}}(r)-U^{\bar{\gamma}_{\bar{t}},x}%
(r)))d\theta,\\
\hat{\alpha}_{\gamma_{t},\bar{\gamma}_{\bar{t}}}(r) & =h(W_{r}^{{\gamma}_{{t}%
}},X^{\bar{\gamma}_{\bar{t}},x}(r),Y^{\bar{\gamma}_{\bar{t}},x}(r),Z^{\bar
{\gamma}_{\bar{t}},x}(r))-h(W_{r}^{\bar{\gamma}_{\bar{t}}},X^{\bar{\gamma
}_{\bar{t}},x}(r),Y^{\bar{\gamma}_{\bar{t}},x}(r),Z^{\bar{\gamma}_{\bar{t}}%
,x}(r)),\\
\hat{\theta}_{\gamma_{t},\bar{\gamma}_{\bar{t}}}(r) & =\int_{0}^{1}%
\frac{\partial h}{\partial x}(W_{r}^{\gamma_{t}},U^{\bar{\gamma}_{\bar{t}}%
,x}(r)+\theta(U^{{\gamma}_{{t}}}(r)-U^{\bar{\gamma}_{\bar{t}},x}%
(r)))d\theta,\\
\hat{\beta}_{\gamma_{t},\bar{\gamma}_{\bar{t}}}(r) & =\int_{0}^{1}%
\frac{\partial h}{\partial y}(W_{r}^{\gamma_{t}},U^{\bar{\gamma}_{\bar{t}}%
,x}(r)+\theta(U^{{\gamma}_{{t}}}(r)-U^{\bar{\gamma}_{\bar{t}},x}%
(r)))d\theta,\\
\hat{\delta}_{\gamma_{t},\bar{\gamma}_{\bar{t}}}(r) & =\int_{0}^{1}%
\frac{\partial h}{\partial z}(W_{r}^{\gamma_{t}},U^{\bar{\gamma}_{\bar{t}}%
,x}(r)+\theta(U^{{\gamma}_{{t}}}(r)-U^{\bar{\gamma}_{\bar{t}},x}(r)))d\theta,
\end{array}
\]
where $U^{\gamma_{t},x}=(X^{\gamma_{t},x},Y^{\gamma_{t},x},Z^{\gamma_{t},x})$.
Under Assumptions $2.2$, $3.2$ and $3,3$, using the similar methods as in
Lemma $3.1$, we get the first three inequalities. For the next three
inequalities, we can write $(\Delta_{h}^{i}X^{\gamma_{t},x},\Delta_{h}%
^{i}Y^{\gamma_{t},x},\Delta_{h}^{i}Z^{\gamma_{t},x})$ as the solution of the
following linearized FBSDE:
\begin{align*}
\Delta_{h}^{i}X^{\gamma_{t},x}(s)  &  =\int_{t}^{s}[\frac{1}{h}\alpha
_{{\gamma_{t},{\gamma}_{{t}}}^{h_{e_{i}}}}(r)+\theta_{\gamma_{t},{{\gamma
}_{{t}}}^{h_{e_{i}}}}(r)\Delta_{h}^{i}X^{\gamma_{t},x}(r)+\beta_{\gamma
_{t},{{\gamma}_{{t}}}^{h_{e_{i}}}}(r)\Delta_{h}^{i}Y^{\gamma_{t},x}%
(r)+\delta_{\gamma_{t},{{\gamma}_{{t}}}^{h_{e_{i}}}}\Delta_{h}^{i}%
Z^{\gamma_{t},x}(r)]dr\\
&  +\int_{t}^{s}[\frac{1}{h}\bar{\alpha}_{\gamma_{t},{{\gamma}_{{t}}%
}^{h_{e_{i}}}}(r)+\bar{\theta}_{\gamma_{t},{{\gamma}_{{t}}}^{h_{e_{i}}}%
}(r)\Delta_{h}^{i}X^{\gamma_{t},x}(r)+\bar{\beta}_{\gamma_{t},{{\gamma}_{{t}}%
}^{h_{e_{i}}}}(r)\Delta_{h}^{i}Y^{\gamma_{t},x}(r)+\bar{\delta}_{\gamma
_{t},{{\gamma}_{{t}}}^{h_{e_{i}}}}\Delta_{h}^{i}Z^{\gamma_{t},x}(r)]dW(r),
\end{align*}%
\begin{align*}
\Delta_{h}^{i}Y^{\gamma_{t},x}(s)  &  =\frac{1}{h}(g(W_{T}^{\gamma
_{t}^{h_{e_{i}}}},X^{\gamma_{t}^{h_{e_{i}}},x}(T))-g(W_{T}^{\gamma_{{t}}%
},X^{\gamma_{{t}},x}(T)))-\int_{s}^{T}[\frac{1}{h}\hat{\alpha}_{{\gamma
_{t},{\gamma}_{{t}}}^{h_{e_{i}}}}(r)+\hat{\theta}_{\gamma_{t},{{\gamma}_{{t}}%
}^{h_{e_{i}}}}(r)\Delta_{h}^{i}X^{\gamma_{t},x}(r)\\
&  +\hat{\beta}_{\gamma_{t},{{\gamma}_{{t}}}^{h_{e_{i}}}}(r)\Delta_{h}%
^{i}Y^{\gamma_{t},x}(r)+\hat{\delta}_{\gamma_{t},{{\gamma}_{{t}}}^{h_{e_{i}}}%
}\Delta_{h}^{i}Z^{\gamma_{t},x}(r)]dr-\int_{s}^{T}\Delta_{h}^{i}Z^{\gamma
_{t},x}(r)dW(r),
\end{align*}
Then the same calculus yields that
\[
E[\sup_{s\in\lbrack t,T]}\mid\Delta_{h}^{i}Y^{\gamma_{t},x}(s)\mid^{2}%
+\sup_{s\in\lbrack t,T]}\mid\Delta_{h}^{i}X^{\gamma_{t},x}(s)\mid^{2}+\mid
\int_{t}^{T}\mid\Delta_{h}^{i}Z^{\gamma_{t},x}\mid^{2}dr\mid]\leq
C_{2}(1+\parallel\gamma_{t}\parallel^{q}+\mid h\mid^{q}).
\]
Notice that%
\[%
\begin{array}
[c]{rl}
& \Delta_{h}^{i}X^{\gamma_{t},x}(s)-\Delta_{\bar{h}}^{i}X^{\bar{\gamma}%
_{\bar{t}},x}(s)\\
= & \int_{t}^{s}[\frac{1}{h}{\alpha}_{\gamma_{t},{\gamma}_{{t}}^{{{h}_{e_{i}}%
}}}(r)-\frac{1}{\bar{h}}{\alpha}_{\bar{\gamma}_{\bar{t}},\bar{\gamma}_{\bar
{t}}^{{\bar{h}_{e_{i}}}}}(r)+{\theta}_{\gamma_{t},{{\gamma}_{{t}}}%
^{{{h}_{e_{i}}}}}(r)\Delta_{{h}}^{i}X^{\gamma_{t},x}(r)-{\theta}_{\bar{\gamma
}_{\bar{t}},{\bar{\gamma}_{\bar{t}}^{{\bar{h}_{e_{i}}}}}}(r)\Delta_{\bar{h}%
}^{i}X^{\bar{\gamma}_{\bar{t}},x}(r)\\
& +{\beta}_{\gamma_{t},{{\gamma}_{{t}}}^{{{h}_{e_{i}}}}}(r)\Delta_{{h}}%
^{i}Y^{\gamma_{t},x}(r)-{\beta}_{\bar{\gamma}_{\bar{t}},{\bar{\gamma}_{\bar
{t}}^{{\bar{h}_{e_{i}}}}}}(r)\Delta_{\bar{h}}^{i}Y^{\bar{\gamma}_{\bar{t}}%
,x}(r)+{\delta}_{\gamma_{t},{{\gamma}_{{t}}}^{{{h}_{e_{i}}}}}\Delta_{{h}}%
^{i}Z^{\gamma_{t},x}(r)-{\delta}_{\bar{\gamma}_{\bar{t}},{\bar{\gamma}%
_{\bar{t}}^{{\bar{h}_{e_{i}}}}}}(r)\Delta_{\bar{h}}^{i}Z^{\bar{\gamma}%
_{\bar{t}},x}(r)]dr\\
& +\int_{t}^{s}[\frac{1}{h}\bar{\alpha}_{\gamma_{t},{\gamma}_{{t}}%
^{{{h}_{e_{i}}}}}(r)-\frac{1}{\bar{h}}\bar{\alpha}_{\bar{\gamma}_{\bar{t}%
},\bar{\gamma}_{\bar{t}}^{{\bar{h}_{e_{i}}}}}(r)+\bar{\theta}_{\gamma
_{t},{{\gamma}_{{t}}}^{{{h}_{e_{i}}}}}(r)\Delta_{{h}}^{i}X^{\gamma_{t}%
,x}(r)-\bar{\theta}_{\bar{\gamma}_{\bar{t}},{\bar{\gamma}_{\bar{t}}^{{\bar
{h}_{e_{i}}}}}}(r)\Delta_{\bar{h}}^{i}X^{\bar{\gamma}_{\bar{t}},x}(r)\\
& +\bar{\beta}_{\gamma_{t},{{\gamma}_{{t}}}^{{{h}_{e_{i}}}}}(r)\Delta_{{h}%
}^{i}Y^{\gamma_{t},x}(r)-\bar{\beta}_{\bar{\gamma}_{\bar{t}},{\bar{\gamma
}_{\bar{t}}^{{\bar{h}_{e_{i}}}}}}(r)\Delta_{\bar{h}}^{i}Y^{\bar{\gamma}%
_{\bar{t}},x}(r)+\bar{\delta}_{\gamma_{t},{{\gamma}_{{t}}}^{{{h}_{e_{i}}}}%
}\Delta_{{h}}^{i}Z^{\gamma_{t},x}(r)-\bar{\delta}_{\bar{\gamma}_{\bar{t}%
},{\bar{\gamma}_{\bar{t}}^{{\bar{h}_{e_{i}}}}}}(r)\Delta_{\bar{h}}^{i}%
Z^{\bar{\gamma}_{\bar{t}},x}(r)]dW(r),
\end{array}
\]
and%
\[%
\begin{array}
[c]{rl}
& \Delta_{h}^{i}Y^{\gamma_{t},x}(s)-\Delta_{\bar{h}}^{i}Y^{\bar{\gamma}%
_{\bar{t}},x}(s)\\
= & \frac{1}{h}(g(W_{T}^{\gamma_{t}^{h_{e_{i}}}},X^{\gamma_{t}^{h_{e_{i}}}%
,x}(T))-g(W_{T}^{\gamma_{{t}}},X^{\gamma_{{t}},x}(T)))-\frac{1}{\bar{h}%
}(g(X_{T}^{\bar{\gamma}_{t}^{{\bar{h}}_{e_{i}}}},X^{\bar{\gamma}_{t}^{{\bar
{h}}_{e_{i}}},x}(T))-g(W_{T}^{\bar{\gamma}_{\bar{t}}},X^{\bar{\gamma}_{\bar
{t}},x}(T)))\\
& -\{\int_{s}^{T}[\frac{1}{h}\hat{\alpha}_{\gamma_{t},{\gamma}_{{t}}%
^{{{h}_{e_{i}}}}}(r)-\frac{1}{\bar{h}}\hat{\alpha}_{\bar{\gamma}_{\bar{t}%
},\bar{\gamma}_{\bar{t}}^{{\bar{h}_{e_{i}}}}}(r)+\hat{\theta}_{\gamma
_{t},{{\gamma}_{{t}}}^{{{h}_{e_{i}}}}}(r)\Delta_{{h}}^{i}X^{\gamma_{t}%
,x}(r)-\hat{\theta}_{\bar{\gamma}_{\bar{t}},{\bar{\gamma}_{\bar{t}}^{{\bar
{h}_{e_{i}}}}}}(r)\Delta_{\bar{h}}^{i}X^{\bar{\gamma}_{\bar{t}},x}(r)\\
& +\hat{\beta}_{\gamma_{t},{{\gamma}_{{t}}}^{{{h}_{e_{i}}}}}(r)\Delta_{{h}%
}^{i}Y^{\gamma_{t},x}(r)-\hat{\beta}_{\bar{\gamma}_{\bar{t}},{\bar{\gamma
}_{\bar{t}}^{{\bar{h}_{e_{i}}}}}}(r)\Delta_{\bar{h}}^{i}Y^{\bar{\gamma}%
_{\bar{t}},x}(r)+\hat{\delta}_{\gamma_{t},{{\gamma}_{{t}}}^{{{h}_{e_{i}}}}%
}\Delta_{{h}}^{i}Z^{\gamma_{t},x}(r)-\hat{\delta}_{\bar{\gamma}_{\bar{t}%
},{\bar{\gamma}_{\bar{t}}^{{\bar{h}_{e_{i}}}}}}(r)\Delta_{\bar{h}}^{i}%
Z^{\bar{\gamma}_{\bar{t}},x}(r)]dr\}\\
& -\int_{s}^{T}(\Delta_{h}^{i}Z^{\gamma_{t},x}(r)-\Delta_{\bar{h}}^{i}%
Z^{\bar{\gamma}_{\bar{t}},x}(r))dW(r).
\end{array}
\]
Then
\[
(\tilde{X}(s),\tilde{Y}(s),\tilde{Z}(s)):=(\Delta_{h}^{i}X^{\gamma_{t}%
,x}(s)-\Delta_{\bar{h}}^{i}X^{\bar{\gamma}_{\bar{t}},x}(s),\Delta_{h}%
^{i}Y^{\gamma_{t},x}(s)-\Delta_{\bar{h}}^{i}Y^{\bar{\gamma}_{\bar{t}}%
,x}(s),\Delta_{h}^{i}Z^{\gamma_{t},x}(s)-\Delta_{\bar{h}}^{i}Z^{\bar{\gamma
}_{\bar{t}},x}(s))
\]
solves the FBSDE
\begin{align*}
\tilde{X}(s)  &  =\int_{t}^{s}[{\theta}_{\gamma_{t},{{\gamma}_{{t}}}%
^{{{h}_{e_{i}}}}}(r)\tilde{X}(r)+{\beta}_{\gamma_{t},{{\gamma}_{{t}}}%
^{{{h}_{e_{i}}}}}(r)\tilde{Y}(r)+{\delta}_{\gamma_{t},{{\gamma}_{{t}}}%
^{{{h}_{e_{i}}}}}\tilde{Z}(r)+\tilde{b}(r)]dr\\
&  +\int_{t}^{s}[\bar{\theta}_{\gamma_{t},{{\gamma}_{{t}}}^{{{h}_{e_{i}}}}%
}(r)\tilde{Y}(r)+\bar{\beta}_{\gamma_{t},{{\gamma}_{{t}}}^{{{h}_{e_{i}}}}%
}(r)\tilde{Y}(r)+\bar{\delta}_{\gamma_{t},{{\gamma}_{{t}}}^{{{h}_{e_{i}}}}%
}\tilde{Z}(r)+\tilde{\sigma}(r)]dW(r),
\end{align*}%
\begin{align*}
\tilde{Y}(s)  &  =\frac{1}{h}(g(W_{T}^{\gamma_{t}^{h_{e_{i}}}},X^{\gamma
_{t}^{h_{e_{i}}},x}(T))-g(W_{T}^{\gamma_{{t}}},X^{\gamma_{{t}},x}%
(T)))-\frac{1}{\bar{h}}(g(W_{T}^{\bar{\gamma}_{t}^{{\bar{h}}_{e_{i}}}}%
,X^{\bar{\gamma}_{t}^{{\bar{h}}_{e_{i}}},x}(T))-g(W_{T}^{\bar{\gamma}_{\bar
{t}}},X^{\bar{\gamma}_{\bar{t}},x}(T)))\\
&  -\int_{s}^{T}[{\hat{\theta}}_{\gamma_{t},{{\gamma}_{{t}}}^{{{h}_{e_{i}}}}%
}(r)\tilde{X}(r)+{\hat{\beta}}_{\gamma_{t},{{\gamma}_{{t}}}^{{{h}_{e_{i}}}}%
}(r)\tilde{Y}(r)+{\hat{\delta}}_{\gamma_{t},{{\gamma}_{{t}}}^{{{h}_{e_{i}}}}%
}\tilde{Z}(r)+\tilde{h}(r)]dr-\int_{s}^{T}\tilde{Z}(r)dW(r),
\end{align*}
where
\begin{align*}
\tilde{b}(r)  &  =[{\beta}_{\gamma_{t},{{\gamma}_{{t}}}^{{{h}_{e_{i}}}}%
}(r)-{\beta}_{\bar{\gamma}_{\bar{t}},{\bar{\gamma}_{\bar{t}}^{{\bar{h}_{e_{i}%
}}}}}(r)]\Delta_{\bar{h}}^{i}Y^{\bar{\gamma}_{\bar{t}},x}(r)+{\theta}%
_{\gamma_{t},{{\gamma}_{{t}}}^{{{h}_{e_{i}}}}}(r)-{\theta}_{\bar{\gamma}%
_{\bar{t}},{\bar{\gamma}_{\bar{t}}^{{\bar{h}_{e_{i}}}}}}(r)]\Delta_{\bar{h}%
}^{i}X^{\bar{\gamma}_{\bar{t}},x}(r)\\
&  +[{\delta}_{\gamma_{t},{{\gamma}_{{t}}}^{{{h}_{e_{i}}}}}(r)-{\delta}%
_{\bar{\gamma}_{\bar{t}},{\bar{\gamma}_{\bar{t}}^{{\bar{h}_{e_{i}}}}}%
}(r)]\Delta_{\bar{h}}^{i}Z^{\bar{\gamma}_{\bar{t}},x}(r)+\frac{1}{h}{\alpha
}_{\gamma_{t},{\gamma}_{{t}}^{{{h}_{e_{i}}}}}(r)-\frac{1}{\bar{h}}{\alpha
}_{\bar{\gamma}_{\bar{t}},\bar{\gamma}_{\bar{t}}^{{\bar{h}_{e_{i}}}}}(r),\\
\tilde{\sigma}(r)  &  =\bar{\beta}_{\gamma_{t},{{\gamma}_{{t}}}^{{{h}_{e_{i}}%
}}}(r)-{\bar{\beta}}_{\bar{\gamma}_{\bar{t}},{\bar{\gamma}_{\bar{t}}^{{\bar
{h}_{e_{i}}}}}}(r)]\Delta_{\bar{h}}^{i}Y^{\bar{\gamma}_{\bar{t}},x}%
(r)+\bar{\beta}_{\gamma_{t},{{\gamma}_{{t}}}^{{{h}_{e_{i}}}}}(r)-{\bar{\beta}%
}_{\bar{\gamma}_{\bar{t}},{\bar{\gamma}_{\bar{t}}^{{\bar{h}_{e_{i}}}}}%
}(r)]\Delta_{\bar{h}}^{i}X^{\bar{\gamma}_{\bar{t}},x}(r)\\
&  +[{\bar{\delta}}_{\gamma_{t},{{\gamma}_{{t}}}^{{{h}_{e_{i}}}}}%
(r)-{\bar{\delta}}_{\bar{\gamma}_{\bar{t}},{\bar{\gamma}_{\bar{t}}^{{\bar
{h}_{e_{i}}}}}}(r)]\Delta_{\bar{h}}^{i}Z^{\bar{\gamma}_{\bar{t}},x}%
(r)\}+\frac{1}{h}{\bar{\alpha}}_{\gamma_{t},{\gamma}_{{t}}^{{{h}_{e_{i}}}}%
}(r)-\frac{1}{\bar{h}}{\bar{\alpha}}_{\bar{\gamma}_{\bar{t}},\bar{\gamma
}_{\bar{t}}^{{\bar{h}_{e_{i}}}}}(r),\\
\tilde{h}(r)  &  =[{\hat{\beta}}_{\gamma_{t},{{\gamma}_{{t}}}^{{{h}_{e_{i}}}}%
}(r)-{\hat{\beta}}_{\bar{\gamma}_{\bar{t}},{\bar{\gamma}_{\bar{t}}^{{\bar
{h}_{e_{i}}}}}}(r)]\Delta_{\bar{h}}^{i}Y^{\bar{\gamma}_{\bar{t}},x}%
(r)+{\hat{\beta}}_{\gamma_{t},{{\gamma}_{{t}}}^{{{h}_{e_{i}}}}}(r)-{\hat
{\beta}}_{\bar{\gamma}_{\bar{t}},{\bar{\gamma}_{\bar{t}}^{{\bar{h}_{e_{i}}}}}%
}(r)]\Delta_{\bar{h}}^{i}X^{\bar{\gamma}_{\bar{t}},x}(r)\\
&  +[{\hat{\delta}}_{\gamma_{t},{{\gamma}_{{t}}}^{{{h}_{e_{i}}}}}%
(r)-{\hat{\delta}}_{\bar{\gamma}_{\bar{t}},{\bar{\gamma}_{\bar{t}}^{{\bar
{h}_{e_{i}}}}}}(r)]\Delta_{\bar{h}}^{i}Z^{\bar{\gamma}_{\bar{t}},x}%
(r)+\frac{1}{h}{\hat{\alpha}}_{\gamma_{t},{\gamma}_{{t}}^{{{h}_{e_{i}}}}%
}(r)-\frac{1}{\bar{h}}{\hat{\alpha}}_{\bar{\gamma}_{\bar{t}},\bar{\gamma
}_{\bar{t}}^{{\bar{h}_{e_{i}}}}}(r).
\end{align*}
Thus, under Assumption 2.2, 3.2 and 3.3, similarly as in Lemma 3.1, we can get
the last three inequalities.\quad$\Box$

\begin{theorem}
For each $\gamma_{t}\in\Lambda$, $\{Y^{{\gamma_{t}^{z},x}}(s),s\in
\lbrack0,T],z\in\mathbb{R}^{d}\}$ has a version which is a.e. of class
$C^{0,2}([0,T]\times\mathbb{R}^{d})$.
\end{theorem}

\noindent\textbf{{Proof.}} In order to simplify the presentation, we only
consider $n=d=1$. Applying Lemma $3.1$, then for each $h,\bar{h}\in
\mathbb{R}\setminus\{0\}$ and $k,\bar{k}\in\mathbb{R}$,\newline%
\[%
\begin{array}
[c]{rl}
& E[\sup_{u\in\lbrack t,T]}\mid Y^{\gamma_{t}^{k},x}(u)-Y^{{\gamma}_{{t}%
}^{\bar{k}},x}(u)\mid^{2}]\leq C_{2}(1+\parallel\gamma_{t}\parallel^{q})\mid
k-\bar{k}\mid^{2},\\
& E[\sup_{u\in\lbrack t,T]}\mid X^{\gamma_{t}^{k},x}(u)-X^{{\gamma}_{{t}%
}^{\bar{k}},x}(u)\mid^{2}]\leq C_{2}(1+\parallel\gamma_{t}\parallel^{q})\mid
k-\bar{k}\mid^{2},\\
& E[\int_{t}^{T}\mid Z^{\gamma_{t}^{k},x}(u)-Z^{{\gamma}_{{t}}^{\bar{k}}%
,x}\mid^{2}du]\leq C_{2}(1+\parallel\gamma_{t}\parallel^{q})\mid k-\bar{k}%
\mid^{2},\\
& E[\sup_{u\in\lbrack t,T]}\mid\Delta_{h}^{i}Y^{\gamma_{t}^{k},x}%
(u)-\Delta_{\bar{h}}^{i}Y^{{\gamma}_{{{t}}}^{\bar{k}},x}(u)\mid^{2}]\\
\leq & C_{2}(1+\parallel\gamma_{t}\parallel^{q}+\parallel\bar{\gamma}_{\bar
{t}}\parallel^{q}+\mid h\mid^{q}+\mid\bar{h}\mid^{q}))(\mid k-\bar{k}\mid
^{2}+\mid h-\bar{h}\mid^{2}),\\
& E[\sup_{u\in\lbrack t,T]}\mid\Delta_{h}^{i}X^{\gamma_{t}^{k},x}%
(u)-\Delta_{\bar{h}}^{i}X^{\bar{\gamma}_{{{t}}}^{\bar{k}},x}(u)\mid^{2}]\\
\leq & C_{2}(1+\parallel\gamma_{t}\parallel^{q}+\parallel\bar{\gamma}_{\bar
{t}}\parallel^{q}+\mid h\mid^{q}+\mid\bar{h}\mid^{q}))(\mid k-\bar{k}\mid
^{2}+\mid h-\bar{h}\mid^{2}),\\
& E[\mid\int_{t}^{T}\mid\Delta_{h}^{i}Z^{\gamma_{t}^{k},x}(u)-\Delta_{\bar{h}%
}^{i}Z^{\gamma_{{{t}}}^{\bar{k}},x}(u)\mid^{2}du\mid]\\
\leq & C_{2}(1+\parallel\gamma_{t}\parallel^{q}+\parallel\bar{\gamma}_{\bar
{t}}\parallel^{q}+\mid h\mid^{q}+\mid\bar{h}\mid^{q}))(\mid k-\bar{k}\mid
^{2}+\mid h-\bar{h}\mid^{2}).
\end{array}
\]
By kolmogorov's criterion, there exists a continuous derivative of
$Y^{\gamma_{t}^{z},x}(s)$ ($X^{\gamma_{t}^{z},x}(s)$) with respect to $z.$
There also exists a mean-square derivative of $Z^{\gamma_{t}^{z},x}(s)$ with
respect to $z$, which is mean square continuous in $z$. We denote them by
\[
(D_{z}Y^{\gamma_{t},x},D_{z}X^{\gamma_{t},x},D_{z}Z^{\gamma_{t},x}).
\]
By Theorem 3.1 and Definition $3.1$, $(D_{z}Y^{\gamma_{t},x},D_{z}%
X^{\gamma_{t},x},D_{z}Z^{\gamma_{t},x})$ is the solution of the following
FBSDE:
\[%
\begin{array}
[c]{cl}
& D_{z}X^{\gamma_{t},x}(s)\\
= & \int_{t}^{s}[b_{\gamma_{t}}^{\prime}(W_{r}^{\gamma_{t}},X^{\gamma_{t}%
,x}(r),Y^{\gamma_{t},x}(r),Z^{\gamma_{t},x}(r))+b_{x}^{\prime}(W_{r}%
^{\gamma_{t}},X^{\gamma_{t},x}(r),Y^{\gamma_{t},x}(r),Z^{\gamma_{t}%
,x}(r))D_{z}X^{\gamma_{t},x}(r)\\
& +b_{y}^{\prime}(W_{r}^{\gamma_{t}},X^{\gamma_{t},x}(r),Y^{\gamma_{t}%
,x}(r),Z^{\gamma_{t},x}(r))D_{z}Y^{\gamma_{t},x}(r)+b_{z}^{\prime}%
(W_{r}^{\gamma_{t}},X^{\gamma_{t},x}(r),Y^{\gamma_{t},x}(r),Z^{\gamma_{t}%
,x}(r))D_{z}Z^{\gamma_{t},x}(r)]dr\\
& +\int_{t}^{s}[\sigma_{\gamma_{t}}^{\prime}(W_{r}^{\gamma_{t}},X^{\gamma
_{t},x}(r),Y^{\gamma_{t},x}(r),Z^{\gamma_{t},x}(r))+\sigma_{x}^{\prime}%
(W_{r}^{\gamma_{t}},X^{\gamma_{t},x}(r),Y^{\gamma_{t},x}(r),Z^{\gamma_{t}%
,x}(r))D_{z}X^{\gamma_{t},x}(r)\\
& +\sigma_{y}^{\prime}(W_{r}^{\gamma_{t}},X^{\gamma_{t},x}(r),Y^{\gamma_{t}%
,x}(r),Z^{\gamma_{t},x}(r))D_{z}Y^{\gamma_{t},x}(r)+\sigma_{z}^{\prime}%
(W_{r}^{\gamma_{t}},X^{\gamma_{t},x}(r),Y^{\gamma_{t},x}(r),Z^{\gamma_{t}%
,x}(r))D_{z}Z^{\gamma_{t},x}(r)]dW(r),
\end{array}
\]%
\[%
\begin{array}
[c]{cl}
& D_{z}Y^{\gamma_{t},x}(s)\\
= & g_{\gamma_{t}}^{\prime}(W_{T}^{\gamma_{t}},X^{\gamma_{t},x}(T))+g_{x}%
^{\prime}(W_{T}^{\gamma_{t}},X^{\gamma_{t},x}(T))D_{z}X^{\gamma_{t},x}%
(T)-\int_{s}^{T}[h_{\gamma_{t}}^{\prime}(W_{r}^{\gamma_{t}},X^{\gamma_{t}%
,x}(r),Y^{\gamma_{t},x}(r),Z^{\gamma_{t},x}(r))\\
& +h_{x}^{\prime}(W_{r}^{\gamma_{t}},X^{\gamma_{t},x}(r),Y^{\gamma_{t}%
,x}(r),Z^{\gamma_{t},x}(r))D_{z}X^{\gamma_{t},x}(r)+h_{y}^{\prime}%
(W_{r}^{\gamma_{t}},X^{\gamma_{t},x}(r),Y^{\gamma_{t},x}(r),Z^{\gamma_{t}%
,x}(r))D_{z}Y^{\gamma_{t},x}(r)\\
& +h_{z}^{\prime}(W_{r}^{\gamma_{t}},X^{\gamma_{t},x}(r),Y^{\gamma_{t}%
,x}(r),Z^{\gamma_{t},x}(r))D_{z}Z^{\gamma_{t},x}(r)]dr-\int_{s}^{T}%
D_{z}Z^{\gamma_{t},x}(r)dW(r).
\end{array}
\]
It is easy to check that the above FBSDE satisfies Assumptions 2.2, 3.2 and
3.3. Then the above FBSDE has a uniqueness solution. Thus, the existence of a
continuous second order derivative of $Y^{\gamma_{t}^{z}}(s)$ with respect to
$z$ can be proved in a similar way.\quad$\Box$

Define
\[
u(\gamma_{t},x):=Y^{\gamma_{t},x}(t),\quad for\quad\gamma_{t}\in\Lambda.
\]
We have the following results about $u(\gamma_{t},x)$.

\begin{lemma}
$\forall t\leq s \leq T$, we have $u(W_{s}^{\gamma_{t}},X^{\gamma_{t}%
,x}(s))=Y^{\gamma_{t},x}(s)$.
\end{lemma}

\noindent\textbf{{Proof.}} For given $\gamma_{t_{1}}$, $t_{1}<t$, set
$X(t_{1})=x$. Consider the solution of FBSDE (2.1) and (2.2) on $[t,T]$:%
\[%
\begin{array}
[c]{rl}%
X^{\gamma_{t_{1}},x}(s)= & X^{\gamma_{t_{1}},x}(t)+\int_{t}^{s}b(W_{r}%
^{\gamma_{t_{1}}},Y^{\gamma_{t_{1}},x}(r),Z^{\gamma_{t_{1}},x}(r))dr+\int%
_{t}^{s}\sigma(W_{r}^{\gamma_{t_{1}}},X^{\gamma_{t_{1}},x}(r),Y^{\gamma
_{t_{1}},x}(r),Z^{\gamma_{t_{1}},x}(r))dW(r),\\
Y^{\gamma_{t_{1}},x}(s)= & g(W_{T}^{\gamma_{t_{1}}},X^{\gamma_{t_{1}}%
,x}(T))-\int_{s}^{T}h(W_{r}^{\gamma_{t_{1}}},X^{\gamma_{t_{1}},x}%
(r),Y^{\gamma_{t_{1}},x}(r),Z^{\gamma_{t_{1}},x}(r))dr-\int_{s}^{T}%
Z^{\gamma_{t_{1}},x}(r)dW(r),\quad s\in\lbrack t,T].
\end{array}
\]
For simplicity, set
\[
\xi=X^{\gamma_{t_{1}},x}(t).
\]
Then we need to prove $u(W_{t}^{\gamma_{t_{1}}},\xi)=Y^{\gamma_{t_{1}},x}(t)$.
Define%
\begin{align*}
W_{t}^{N,\gamma_{t_{1}}}  &  :=\Sigma_{i=1}^{N}I_{A_{i}}x_{t}^{i},\\
\xi^{N}  &  :=\Sigma_{i=1}^{N}I_{A_{i}}a^{i},
\end{align*}
where $\{A_{i}\}_{i=1}^{N}$ is a division of $\mathcal{F}_{t}$, $x_{t}^{i}\in
A_{i}\cap{\Lambda}$, $i=1,2,\cdots,N$. For any $i$, $(X^{x_{t}^{i},a^{i}%
}(s),Y^{x_{t}^{i},a^{i}}(s),Y^{x_{t}^{i},a^{i}}(s))$ is the solution of the
following FBSDE:%
\[%
\begin{array}
[c]{rl}%
X^{x_{t}^{i},a^{i}}(s)= & a^{i}+\int_{t}^{s}b({x_{t}^{i}},X^{x_{t}^{i},a^{i}%
}(r),Y^{x_{t}^{i},a^{i}}(r),Z^{x_{t}^{i}}(r))dr+\int_{t}^{s}\sigma({x_{t}^{i}%
},X^{x_{t}^{i},a^{i}}(r),Y^{x_{t}^{i},a^{i}}(r),Z^{x_{t}^{i},a^{i}%
}(r))dW(r),\\
Y^{x_{t}^{i},a^{i}}(s)= & g(W_{T}^{x_{t}^{i}},X^{x_{t}^{i},a^{i}}(T))-\int%
_{s}^{T}h({x_{t}^{i}},X^{x_{t}^{i},a^{i}}(r),Y^{x_{t}^{i},a^{i}}%
(r),Z^{x_{t}^{i},a^{i}}(r))dr-\int_{s}^{T}Z^{x_{t}^{i},a^{i}}(r)dW(r),\quad
s\in\lbrack t,T].
\end{array}
\]
Multiplying by $I_{A_{i}}$ and adding the corresponding terms, we obtain:
\begin{align*}
\sum_{i=1}^{N}I_{A_{i}}X^{x_{t}^{i},a^{i}}(s)  &  =\sum_{i=1}^{N}I_{A_{i}%
}a^{i}+\int_{t}^{s}b(\sum_{i=1}^{N}I_{A_{i}}{x_{t}^{i}},\sum_{i=1}^{N}%
I_{A_{i}}Y^{x_{t}^{i},a^{i}}(r),\sum_{i=1}^{N}I_{A_{i}}Y^{x_{t}^{i},a^{i}%
}(r),\sum_{i=1}^{N}I_{A_{i}}Z^{x_{t}^{i},a^{i}}(r))dr\\
&  +\int_{t}^{s}\sigma(\sum_{i=1}^{N}I_{A_{i}}{x_{t}^{i}},\sum_{i=1}%
^{N}I_{A_{i}}Y^{x_{t}^{i},a^{i}}(r),\sum_{i=1}^{N}I_{A_{i}}Y^{x_{t}^{i},a^{i}%
}(r),\sum_{i=1}^{N}I_{A_{i}}Z^{x_{t}^{i},a^{i}}(r))dW(r),
\end{align*}%
\begin{align*}
\sum_{i=1}^{N}I_{A_{i}}Y^{x_{t}^{i},a^{i}}(s)  &  =g(\sum_{i=1}^{N}I_{A_{i}%
}W_{T}^{x_{t}^{i}},\sum_{i=1}^{N}I_{A_{i}}X^{x_{t}^{i},a^{i}}(T))-\int_{s}%
^{T}\sum_{i=1}^{N}I_{A_{i}}Z^{x_{t}^{i},a^{i}}(r)dW(r),\\
&  -\int_{s}^{T}h(\sum_{i=1}^{N}I_{A_{i}}{x_{t}^{i}},\sum_{i=1}^{N}I_{A_{i}%
}X^{x_{t}^{i},a^{i}}(r),\sum_{i=1}^{N}I_{A_{i}}Y^{x_{t}^{i},a^{i}}%
(r),\sum_{i=1}^{N}I_{A_{i}}Z^{x_{t}^{i},a^{i}}(r))dr.\quad s\in\lbrack t,T].
\end{align*}
So
\[
X^{\gamma_{t},x}(s)=\sum_{i=1}^{N}I_{A_{i}}X^{x_{t}^{i},a^{i}}(s)\text{,\ }%
Y^{\gamma_{t},x}(s)=\sum_{i=1}^{N}I_{A_{i}}Y^{x_{t}^{i},a^{i}}(s),Z^{\gamma
_{t},a^{i}}(s)=\sum_{i=1}^{N}I_{A_{i}}Z^{x_{t}^{i},a^{i}}(s)
\]
is the solution of the above FBSDE. By the definition of $u$, we get
\[
Y^{^{W_{t}^{N,\gamma_{t_{1}}}},\xi^{N}}(t)=\sum_{i=1}^{N}I_{A_{i}}Y^{x_{t}%
^{i},a^{i}}(t)=\sum_{i=1}^{N}I_{A_{i}}u(x_{t}^{i},a^{i})=u(W_{t}%
^{N,\gamma_{t_{1}}},\xi^{N}).
\]
For the general case, following the method in Peng and Wang \cite{Peng S 3}
(Lemma $4.3$), we choose a simple adapted process $\{\gamma_{t}^{i}%
\}_{i=1}^{\infty}$ such that $E\parallel\gamma_{t}^{i}-W_{t}^{\gamma_{t_{1}}%
}\parallel$ convergence to $0$ as $i\rightarrow\infty$. Using the same
procedure for $\xi$, we obtain\newline%
\[
E\mid Y^{\gamma_{t}^{i},\xi^{i}}(t)-Y^{W_{t}^{\gamma_{t_{1}}},\xi}(t)\mid
^{2}\leq CE[\parallel\gamma_{t}^{i}-W_{t}^{\gamma_{t_{1}}}\parallel+\mid
\xi^{i}-\xi\mid]
\]
and
\[
E\mid u(\gamma_{t}^{i},\xi^{i})-u(W_{t}^{\gamma_{t_{1}}},\xi)\mid^{2}\leq
\bar{C}E[\parallel\gamma_{t}^{i}-W_{t}^{\gamma_{t_{1}}}\parallel+\mid\xi
^{i}-\xi\mid].
\]
This completes the proof.\quad$\Box$

By Theorem 3.1, 3.2 and the definition of Dupire's vertical derivative, we
have the following corollary.

\begin{corollary}
$D_{z}u(\gamma_{t},x)$ and $D_{zz}u(\gamma_{t},x)$ exist. Moreover,
$u(\gamma_{t},x)$, $D_{z}u(\gamma_{t},x)$ and $D_{zz}u(\gamma_{t},x)$ are
$\Lambda$-continuous.
\end{corollary}

\noindent\textbf{{Proof.}} By Theorem $3.2$ we know that $D_{z}u(\gamma
_{t},x)$ and $D_{zz}u(\gamma_{t},x)$ exist. In the following, we only prove
$u(\gamma_{t},x)$ is $\Lambda$-continuous. The proof for the continuous
property of $D_{z}u(\gamma_{t},x)$ and $D_{zz}u(\gamma_{t},x)$ is similar.
Taking expectation on both sides of equation (2.2),
\[
u({\gamma_{t},x})=Eg(W_{T}^{\gamma_{t}},X^{\gamma_{t},x}(T))-E\int_{t}%
^{T}h(W_{r}^{\gamma_{t}},X^{\gamma_{t},x}(r),Y^{\gamma_{t},x}(r),Z^{\gamma
_{t},x}(r))dr.
\]
For $\gamma_{t},\bar{\gamma}_{\bar{t}}\in\Lambda,$ $\bar{t}\geq t$, we have
\[%
\begin{array}
[c]{rl}
& \mid u(\gamma_{t},x)-u(\bar{\gamma}_{\bar{t}},x)\mid\\
\leq & E[\mid g(W_{T}^{\gamma_{t}},X^{\gamma_{t},x}(T))-g(W_{T}^{\bar{\gamma
}_{\bar{t}}},X^{\bar{\gamma}_{\bar{t}},x}(T))\mid]+E[\int_{t}^{\bar{t}}\mid
h(W_{r}^{{\gamma}_{{t}}},X^{{\gamma}_{{t}},x}(r),Y^{{\gamma}_{{t}}%
,x}(r),Z^{{\gamma}_{{t}},x}(r))\mid dr]\\
& +E[\int_{\bar{t}}^{T}\mid h(W_{r}^{\gamma_{t}},X^{\gamma_{t}}(r),Y^{\gamma
_{t},x}(r),Z^{\gamma_{t},x}(r))-h(W_{r}^{\bar{\gamma}_{\bar{t}}}%
,X^{\bar{\gamma}_{\bar{t}},x}(r),Y^{\bar{\gamma}_{\bar{t}},x}(r),Z^{\bar
{\gamma}_{\bar{t}},x}(r))\mid dr]\\
\leq & E[C_{1}(1+\parallel W_{T}^{\gamma_{t}}\parallel^{k}+\parallel
W_{T}^{\bar{\gamma}_{\bar{t}}}\parallel^{k})\parallel{\gamma_{t}}-{\bar
{\gamma}_{\bar{t}}}\parallel\\
& +3(\bar{t}-t)^{\frac{1}{2}}(\int_{t}^{\bar{t}}(\mid h(W^{{\gamma}_{{t}}%
}(r),0,0,0)\mid^{2}+\mid CX^{\gamma_{t},x}(r)\mid^{2}+\mid CZ^{\gamma_{t}%
,x}(r)\mid^{2}+\mid CY^{\gamma_{t},x}(r)\mid^{2})dr)^{\frac{1}{2}}\\
& +C\int_{\bar{t}}^{T}(\mid X^{\gamma_{t},x}(r)-X^{\bar{\gamma}_{\bar{t}}%
,x}(r)\mid+\mid Y^{\gamma_{t},x}(r)-Y^{\bar{\gamma}_{\bar{t}},x}(r)\mid+\mid
Z^{\gamma_{t},x}(r)-Z^{\bar{\gamma}_{\bar{t}},x}(r)\mid)dr].
\end{array}
\]
By Theorem 3.1, for some constant $C_{1}$ depending only in $C,k,x$ and $T$,
\[
\mid u(\gamma_{t},x)-u(\bar{\gamma}_{\bar{t}},x)\mid\leq C_{1}(1+\parallel
\gamma_{t}\parallel^{k}+\parallel\bar{\gamma}_{\bar{t}}\parallel
^{k})(\parallel\gamma_{t}-\bar{\gamma}_{\bar{t}}\parallel+\mid t-\bar{t}%
\mid^{\frac{1}{2}}).
\]
This completes the proof.\quad$\Box$

Using similar methods in this section, we can also prove that there exists the
second order partial derivative of $u$ in $x$. Finally, we have $u\in
\mathbb{C}^{0,2,2}(\Lambda\times\mathbb{R}^{n})$.

\subsection{Path regularity of process Z}

In Pardoux and Peng \cite{Pardoux.E 2}, BSDE is only state-dependent, i.e.,
$h=h(t,\gamma(t),y,z)$ and $g=g(\gamma(T))$. Under appropriate assumptions,
$Y$ and $Z$ are related in the following sense:
\[
Z^{\gamma_{t}}(s)=\nabla_{x}u(s,\gamma_{t}(t)+W(s)-W(t)),\quad P-a.s.
\]
Peng and Wang \cite{Peng S 3} extends this result to the path-dependent case.
The corresponding BSDE is
\[
Y^{\gamma_{t}}(s)=g(W_{T}^{\gamma_{t}})-\int_{s}^{T}h(W_{r}^{\gamma_{t}%
},Y^{\gamma_{t}}(r),Z^{\gamma_{t}}(r))dr-\int_{s}^{T}Z^{\gamma_{t}%
}(r)dW(r),\quad s\in\lbrack t,T].
\]
where $W_{T}^{\gamma_{t}}=I_{s\leq t}\gamma_{t}(s)+I_{t<s\leq T}(\gamma
_{t}(t)+W(s)-W(t)).$ Then under some assumptions, they obtained
\[
Z^{\gamma_{t},x}(s)=D_{z}u(W_{s}^{\gamma_{t}}),\quad P-a.s.
\]
In this paper, we generalize it to the Non-Markovian fully coupled FBSDE case.

\begin{theorem}
Under Assumptions $2.2,$ $3.2$ and $3.3$, for each $\gamma_{t}\in\Lambda$, the
process $(Z^{\gamma_{t},x}(s))_{s\in\lbrack t,T]}$ has a continuous version
with the form,
\[
Z^{\gamma_{t},x}(s)=\sigma(W_{s}^{\gamma_{t}},X^{\gamma_{t},x}(s),Y^{\gamma
_{t},x}(s),Z^{\gamma_{t},x}(s))\nabla_{x}u(W_{s}^{\gamma_{t}},X^{\gamma_{t}%
,x}(s))+D_{z}u(W_{s}^{\gamma_{t}},X^{\gamma_{t},x}(s)),\quad P-a.s.
\]

\end{theorem}

To prove the above Theorem, we need the following lemma essentially from
Pardoux and Peng \cite{Pardoux.E 2}.

\begin{lemma}
Let $\gamma_{t}$ and some $\bar{t}\in\lbrack t,T]$ be given. Suppose that
\[
g(\gamma,z)=\varphi(\gamma(\bar{t}),\gamma(T)-\gamma(\bar{t}),z),
\]
where $\varphi$ is in $C_{p}^{3}(\mathbb{R}^{2d}\times\mathbb{R}%
^{m};\mathbb{R}^{m})$. For $\phi=b,\sigma$ and $h$, suppose that
\[
\phi(\gamma_{t},x,y,z)=\bar{\phi}_{1}(s,\gamma_{s}(s),x,y,z)I_{[0,\bar{t}%
)}(s)+\bar{\phi}_{2}(s,\gamma_{s}(\bar{t}),\gamma_{s}(s)-\gamma_{s}(\bar
{t}),x,y,z)I_{[\bar{t},T]}(s),
\]
where $\bar{\phi}_{1},\bar{\phi}_{2}\in C^{0,3}$. Then for each $s\in\lbrack
t,T]$,
\[
Z^{\gamma_{t},x}(s)=\sigma(W_{s}^{\gamma_{t}},X^{\gamma_{t},x}(s),Y^{\gamma
_{t},x}(s),Z^{\gamma_{t},x}(s))\nabla_{x}u(W_{s}^{\gamma_{t}},X^{\gamma_{t}%
,x}(s))+D_{z}u(W_{s}^{\gamma_{t}},X^{\gamma_{t},x}(s)),\quad P-a.s.
\]

\end{lemma}

\noindent\textbf{{Proof.}} We only consider one dimentional case. For
$s\in\lbrack\bar{t},T]$, the FBSDE (2.1) and (2.2) can be rewritten as%
\[%
\begin{array}
[c]{rl}%
X^{\gamma_{s},x}(u)= & x+\int_{{s}}^{u}\bar{b}_{2}(r,\gamma_{s}(\bar
{t}),W^{\gamma_{s}}(r)-\gamma_{s}(\bar{t}),X^{\gamma_{s},x}(r),Y^{\gamma
_{s},x}(r),Z^{\gamma_{s},x}(r))dr\\
& +\int_{{s}}^{u}\bar{\sigma}_{2}(r,\gamma_{s}(\bar{t}),W^{\gamma_{s}%
}(r)-\gamma_{s}(\bar{t}),X^{\gamma_{s},x}(r),Y^{\gamma_{s},x}(r),Z^{\gamma
_{s},x}(r))dW(r),\quad u\in\lbrack s,T],\\
Y^{\gamma_{t},x}(u)= & \varphi(\gamma_{s}(\bar{t}),W^{\gamma_{s}}({T}%
)-\gamma_{s}(\bar{t}),X^{\gamma_{s},x}(T))\\
& -\int_{u}^{T}\bar{h}_{2}(r,\gamma_{s}(\bar{t}),W^{\gamma_{s}}(r)-\gamma
_{s}(\bar{t}),X^{\gamma_{s},x}(r),Y^{\gamma_{s},x}(r),Z^{\gamma_{s}%
,x}(r))dr-\int_{u}^{T}Z^{\gamma_{s},x}(r)dW(r),\quad u\in\lbrack s,T].
\end{array}
\]
For $s\in\lbrack t,\bar{t}]$,%
\[%
\begin{array}
[c]{rl}%
X^{\gamma_{s},x}(u)= & x+\int_{{s}}^{u}\bar{b}_{1}(r,W^{\gamma_{s}}%
({r}),X^{\gamma_{s},x}(r),Y^{\gamma_{s},x}(r),Z^{\gamma_{s},x}(r))dr\\
& +\int_{{s}}^{u}\bar{\sigma}_{1}(r,W^{\gamma_{s}}(r),X^{\gamma_{s}%
,x}(r),Y^{\gamma_{s},x}(r),Z^{\gamma_{s},x}(r))dW(r),\quad u\in\lbrack
s,\bar{t}],\\
X^{\gamma_{s},x}(u)= & X^{\gamma_{s},x}(\bar{t})+\int_{{\bar{t}}}^{u}\bar
{b}_{2}(r,W^{\gamma_{s}}(\bar{t}),W^{\gamma_{s}}(r)-W^{\gamma_{s}}(\bar
{t})),X^{\gamma_{s},x}(r),Y^{\gamma_{s},x}(r),Z^{\gamma_{s},x}(r))dr\\
& +\int_{{\bar{t}}}^{u}\bar{\sigma}_{2}(r,W^{\gamma_{s}}(\bar{t}%
),W^{\gamma_{s}}(r)-W^{\gamma_{s}}(\bar{t})),X^{\gamma_{s},x}(r),Y^{\gamma
_{s},x}(r),Z^{\gamma_{s},x}(r))dW(r),\quad u\in\lbrack\bar{t},T],
\end{array}
\]
and%
\[%
\begin{array}
[c]{rl}%
Y^{\gamma_{s},x}(u)= & \varphi(W^{\gamma_{s}}(\bar{t}),W^{\gamma_{s}}%
({T})-W^{\gamma_{s}}(\bar{t}),X^{\gamma_{s},x}(T))\\
& -\int_{u}^{T}\bar{h}_{2}(r,W^{\gamma_{s}}(\bar{t}),W^{\gamma_{s}%
}(r)-W^{\gamma_{s}}(\bar{t}),X^{\gamma_{s},x}(r),Y^{\gamma_{s},x}%
(r),Z^{\gamma_{s},x}(r))dr-\int_{u}^{T}Z^{\gamma_{s},x}(r)dW(r),\quad
u\in\lbrack\bar{t},T],\\
Y^{\gamma_{s},x}(u)= & Y^{\gamma_{s},x}(\bar{t})-\int_{u}^{\bar{t}}\bar{h}%
_{1}(r,W^{\gamma_{s}}(\bar{t}),X^{\gamma_{s},x}(r),Y^{\gamma_{s}%
,x}(r),Z^{\gamma_{s},x}(r))dr-\int_{u}^{\bar{t}}Z^{\gamma_{s},x}(r)dW(r),\quad
u\in\lbrack s,\bar{t}].
\end{array}
\]
Now consider the following quasilinear parabolic differential equations, which
is defined on $[\bar{t},T]\times\mathbb{R}^{2}$ and parameterized by
$x\in\mathbb{R}$,
\[%
\begin{array}
[c]{l}%
\nabla_{s}u_{2}(s,x,y,z)+\mathcal{L}u_{2}(s,x,y,z)+\nabla_{z}\nabla_{y}%
u_{2}\bar{\sigma}_{2}(s,x,y,z,u_{2},v_{2})+\frac{1}{2}\nabla_{yy}u_{2}\\
=\bar{h}_{2}(s,x,y,z,u_{2},\nabla_{y}u_{2}(s,x,y)\bar{\sigma}_{2}%
(s,x,y,z,u_{2},v_{2})),\\
v_{2}(s,x,y,z))=\nabla_{z}u_{2}(s,x,y,z)\bar{\sigma}_{2}(s,x,y,z,u_{2}%
,v_{2})+\nabla_{y}u_{2},\\
u_{2}(T,x,y,z)=\varphi(x,y,z).
\end{array}
\]

where $\mathcal{L}=\frac{1}{2}\bar{\sigma}_{2}^{2}{\nabla_{zz}}+\bar{b}%
_{2}{\nabla_{z}}$. The other one is defined on $[t,\bar{t}]\times\mathbb{R}$,
\[%
\begin{array}
[c]{l}%
\nabla_{s}u_{1}(s,x,z)+\mathcal{L}u_{1}(s,x,z)+\nabla_{z}\nabla_{x}u_{1}%
\bar{\sigma}_{1}(s,x,z,u_{1},v_{1})+\frac{1}{2}\nabla_{xx}u_{1}\\
=\bar{h}_{1}(s,x,z,u_{1},\nabla_{z}u_{1}\bar{\sigma}_{1}(s,x,z,u_{1}%
,v_{1})),\\
v_{1}(s,x,z)=\nabla_{z}u_{1}\bar{\sigma}_{1}(s,x,z,u_{1},v_{1})+\nabla
_{x}u_{1},\\
u_{1}(\bar{t},x,z)=u_{2}(\bar{t},x,0,z).
\end{array}
\]
where $\mathcal{L}=\frac{1}{2}\bar{\sigma}_{1}^{2}{\nabla_{zz}}+\bar{b}%
_{1}{\nabla_{z}}$. Following Theorem $3.1$ , $3.2$ of Paroux-Peng
\cite{Pardoux.E 2}, we have $u_{2}\in C^{1,2}([\bar{t},T]\times\mathbb{R}%
^{2};\mathbb{R})$, $u_{1}\in C^{1,2}([t,\bar{t}]\times\mathbb{R};\mathbb{R})$
and
\[
u(\gamma_{s},z)=u_{1}(s,\gamma_{s}(s),z)I_{[t,\bar{t})}(s)+u_{2}(s,\gamma
_{s}(\bar{t}),\gamma_{s}(s)-\gamma_{s}(\bar{t}),z)I_{[\bar{t},T]}(s).
\]
By It\^{o} formula and the uniqueness theorem of BSDE,%
\[%
\begin{array}
[c]{rl}%
Y^{\gamma_{t},x}(s)= & u_{1}(s,W^{\gamma_{t}}(s),X^{\gamma_{t},x}(s)),\quad
t\leq s<\bar{t},\\
Y^{\gamma_{t},x}(s)= & u_{2}(s,W^{\gamma_{t}}(\bar{t}),W^{\gamma_{t}%
}(s)-W^{\gamma_{t}}(\bar{t}),X^{\gamma_{t},x}(s)),\quad\bar{t}\leq s\leq T,\\
Z^{\gamma_{t},x}(s)= & \nabla_{z}u_{1}(s,W^{\gamma_{t}}(s))\bar{\sigma}%
_{1}(s,W^{\gamma_{t}}(s),X^{\gamma_{t},x}(s),u_{1},Z^{\gamma_{t},x}%
(s))+\nabla_{x}u_{1},\quad t\leq s<\bar{t},\\
Z^{\gamma_{t},x}(s)= & \nabla_{z}u_{2}(s,W^{\gamma_{t}}(\bar{t}),W^{\gamma
_{t}}(s)-W^{\gamma_{t}}(\bar{t}),X^{\gamma_{t},x}(s))\\
& \cdot\bar{\sigma}_{2}(s,W^{\gamma_{t}}(s),W^{\gamma_{s}}(s)-W^{\gamma_{s}%
}(\bar{t})),X^{\gamma_{t},x}(s),u_{2},Z^{\gamma_{t},x}(s))+\nabla_{x}%
u_{2},\quad\bar{t}\leq s\leq T.
\end{array}
\]
Finally, for each $s\in\lbrack t,T]$,
\[
Z^{\gamma_{t},x}(s)=\sigma(W_{s}^{\gamma_{t}},X^{\gamma_{t},x}(s),Y^{\gamma
_{t},x}(s),Z^{\gamma_{t},x}(s))\nabla_{x}u(W_{s}^{\gamma_{t}},X^{\gamma_{t}%
,x}(s))+D_{z}u(W_{s}^{\gamma_{t}},X^{\gamma_{t},x}(s)),\quad P-a.s.
\]
In particular, we have
\[
Z^{\gamma_{t},x}(t)=\sigma(\gamma_{t},x,u(\gamma_{t},x),Z^{\gamma_{t}%
,x}(t))\nabla_{x}u(\gamma_{t},x)+D_{z}u(\gamma_{t},x).\quad\gamma_{t}%
\in\Lambda.
\]
This completes the proof.\quad$\Box$

Now we give the proof of Theorem 3.3.

\noindent\textbf{{Proof}}$\mathbf{{.}}$ For each fixed $t\in\lbrack0,T]$ and
positive integer n, we introduce a mapping $\gamma^{n}(\bar{\gamma}%
_{s},x):\Lambda_{s}\mapsto\Lambda_{s},$
\[
\gamma^{n}(\bar{\gamma}_{s})(r)=\bar{\gamma}_{s}(r)I_{[0,t)}+\Sigma
_{k=0}^{n-1}\bar{\gamma}_{s}(t_{k+1}^{n}\wedge s)I_{[t_{k}^{n}\wedge
s]}(r)+\bar{\gamma}_{s}(s)I_{\{s\}}(r),\;s\in\lbrack0,T],
\]
where $t_{k}^{n}=t+\frac{k(T-t)}{n}$, $k=0,1,\ldots,n.$ Note that $\Psi$
represents $b$, $\sigma$ or $h$. Define
\[
g^{n}(\bar{\gamma},x):=g(\gamma^{n}(\bar{\gamma})),\quad\Psi^{n}(\bar{\gamma
}_{s},x,y,z):=\Psi(\gamma^{n}(\bar{\gamma}_{s},x),x,y,z).
\]
For each $n$, there exists some functions $\varphi_{n}$ defined on
$\Lambda_{t}\times\mathbb{R}^{n\times d}$ and $\psi_{n}$ defined on
$[t,T]\times{\Lambda_{t}}\times\mathbb{R}^{n\times d}\times{\mathbb{R}%
^{m}\times\mathbb{R}^{m\times d}}$ such that%
\[%
\begin{array}
[c]{rl}%
g^{n}(\bar{\gamma},x)= & \varphi_{n}(\bar{\gamma}_{t},\bar{\gamma}(t_{1}%
^{n})-\bar{\gamma}(t),\cdots,\bar{\gamma}(t_{n}^{n})-\bar{\gamma}(t_{n-1}%
^{n}),x),\\
\Psi^{n}(\bar{\gamma}_{s},x,y,z)= & \psi_{n}(s,\bar{\gamma}_{t},\bar{\gamma
}_{s}(t_{1}^{n}\wedge s)-\bar{\gamma}_{s}(t),\cdots,\bar{\gamma}_{s}(t_{n}%
^{n}\wedge s)-\bar{\gamma}_{s}(t_{n-1}^{n}\wedge s),x,y,z).
\end{array}
\]
Indeed, we can set
\begin{align*}
&  \bar{\varphi}_{n}(\bar{\gamma}_{t},x_{1},\cdots,x_{n},x):=g((\bar{\gamma
}_{t}(s)I_{[0,t)}(s)+\Sigma_{k=1}^{n}x_{k}I_{[t_{k-1}^{n},t_{k}^{n})}%
(s)+x_{n}I_{\{T\}}(s))_{0\leq s\leq T},x),\\
&  \varphi_{n}(\bar{\gamma}_{t},x_{1},\cdots,x_{n},x):=\bar{\varphi}_{n}%
(\bar{\gamma}_{t},\bar{\gamma}_{t}+x_{1},\bar{\gamma}_{t}(t)+x_{1}%
+x_{2},\cdots,\bar{\gamma}_{t}(t)+\Sigma_{i=1}^{n}x_{i},x).
\end{align*}
Then by Assumptions $3.2$ and $3.3$, we obtain that $\varphi_{n}(\bar{\gamma
}_{t},x_{1},\cdots,x_{n},x)$ is a $C_{p}^{3}$-function of $x_{1},\cdots,x_{n}$
for each fixed $\bar{\gamma}_{t}$. In particular, for each $\bar{\gamma}%
\in\Lambda,$
\[
\nabla_{x_{i}}\varphi_{n}(\bar{\gamma}_{t},\bar{\gamma}(t_{1}^{n})-\bar
{\gamma}(t),\cdots,\bar{\gamma}(t_{n}^{n})-\bar{\gamma}(t_{n-1}^{n}%
),x)=g_{\gamma_{t_{i-1}^{n}}}^{\prime}(\gamma^{n}(\bar{\gamma}),x).
\]
\newline For any $\bar{t}\geq t$, $\bar{\gamma}_{\bar{t}}\in\Lambda_{\bar{t}}%
$, consider the following BSDEs:
\[
Y^{n,\bar{\gamma}_{\bar{t}}}(s)=g^{n}(W_{T}^{\bar{\gamma}_{\bar{t}}}%
,X^{n,\bar{\gamma}_{\bar{t}}}({T}))-\int_{s}^{T}h^{n}(W_{r}^{\bar{\gamma
}_{\bar{t}}},X^{n,\bar{\gamma}_{\bar{t}}}(r),Y^{n,\bar{\gamma}_{\bar{t}}%
}(r),Z^{n,\bar{\gamma}_{\bar{t}}}(r))dr-\int_{s}^{T}Y^{n,\bar{\gamma}_{\bar
{t}}}(r)dW(r).
\]
we denote
\[
u^{n}(\bar{\gamma}_{\bar{t}},x):=Y^{n,\bar{\gamma}_{\bar{t}}}(t),\quad
\bar{\gamma}_{\bar{t}}\in\Lambda.
\]
Following the argument as in Lemma 3.3, for each $s\in\lbrack t,T]$, we
obtain
\[
Z^{n,\bar{\gamma}_{\bar{t}}}(s)=\sigma^{n}(W_{s}^{\bar{\gamma}_{\bar{t}}%
},X^{n,\bar{\gamma}_{\bar{t}}}(s),Y^{n,\bar{\gamma}_{\bar{t}}}(s),Z^{n,\bar
{\gamma}_{\bar{t}}}(s))\nabla_{x}u^{n}(W_{s}^{\bar{\gamma}_{\bar{t}}%
},X^{n,\bar{\gamma}_{\bar{t}}}(s))+D_{z}u^{n}(W_{s}^{\bar{\gamma}_{\bar{t}}%
},X^{n,\bar{\gamma}_{\bar{t}}}(s)),\quad P-a.s.
\]
Let $C_{0}$ be a constant depends only on $C,T$ and $k,x$, which allowed to
change from line by line. Following the similar calculus as in Lemma $3.1$ and
Theorem $3.1,$ we get%
\[%
\begin{array}
[c]{rl}
& \mid u^{n}(\bar{\gamma}_{\bar{t}},x)-u(\bar{\gamma}_{\bar{t}},x)\mid\\
\leq & C_{0}E[\mid g^{n}(W_{T}^{\bar{\gamma}_{\bar{t}}},X^{\bar{\gamma}%
_{\bar{t}},x}(T))-g(W_{T}^{\bar{\gamma}_{\bar{t}}},X^{\bar{\gamma}_{\bar{t}%
},x}(T))\mid\\
& +\int_{\bar{t}}^{T}\mid h^{n}(W_{r}^{\bar{\gamma}_{\bar{t}}},X^{n,\bar
{\gamma}_{\bar{t}},x}(r),Y^{n,\bar{\gamma}_{\bar{t}},x}(r),Z^{n,\bar{\gamma
}_{\bar{t}},x}(r))-h(W_{r}^{\bar{\gamma}_{\bar{t}}},X^{\bar{\gamma}_{\bar{t}%
},x}(r),Y^{\bar{\gamma}_{\bar{t}},x}(r),Z^{\bar{\gamma}_{\bar{t}},x}%
(r))\mid^{2}dr]^{\frac{1}{2}}\\
\leq & C_{0}(1+\parallel\bar{\gamma}_{\bar{t}}\parallel^{k})(\frac{1}%
{n^{\frac{1}{4}}}+\parallel\gamma^{n}(\bar{\gamma}_{\bar{t}})-\bar{\gamma
}_{\bar{t}}\parallel).
\end{array}
\]
and
\begin{align*}
&  \mid{D_{z}}u^{n}(\bar{\gamma}_{\bar{t}},x)-{D_{z}}u(\bar{\gamma}_{\bar{t}%
},x)\mid\leq C_{0}(1+\parallel\bar{\gamma}_{\bar{t}}\parallel^{k})(\frac
{1}{n^{\frac{1}{4}}}+\parallel\gamma^{n}(\bar{\gamma}_{\bar{t}})-\bar{\gamma
}_{\bar{t}}\parallel),\\
&  \mid{D_{zz}}u^{n}(\bar{\gamma}_{\bar{t}},x)-{D_{zz}}u(\bar{\gamma}_{\bar
{t}},x)\mid\leq C_{0}(1+\parallel\bar{\gamma}_{\bar{t}}\parallel^{k})(\frac
{1}{n^{\frac{1}{4}}}+\parallel\gamma^{n}(\bar{\gamma}_{\bar{t}})-\bar{\gamma
}_{\bar{t}}\parallel).
\end{align*}
Since%
\[%
\begin{array}
[c]{rl}
& \lim_{n}E[\sup_{s\in\lbrack t,T]}\mid D_{z}u^{n}(W_{s}^{\gamma_{t}%
},X^{n,\gamma_{t},x}(s))-D_{z}u(W_{s}^{\gamma_{t}},X^{\gamma_{t},x}%
(s))\mid^{2}]\\
\leq & C_{0}\lim_{n}E[\sup_{s\in\lbrack t,T]}\mid(1+\parallel W_{s}%
^{\gamma_{t}}\parallel^{k}+\parallel X^{n,\gamma_{t},x}(s)\parallel^{k}%
)(\frac{1}{n^{\frac{1}{4}}}+\parallel\gamma^{n}(W_{s}^{\gamma_{t}}%
)-W_{s}^{\gamma_{t}}\parallel)\mid^{2}]\\
\leq & C_{0}\lim_{n}(1+\parallel\bar{\gamma}_{\bar{t}}\parallel^{k})(\frac
{1}{n^{\frac{1}{4}}}+\parallel\gamma^{n}(\bar{\gamma}_{\bar{t}})-\bar{\gamma
}_{\bar{t}}\parallel)=0
\end{array}
\]
and $\lim_{n}E[\mid\int_{t}^{T}\mid Z^{\gamma_{t},x}(u)-Z^{n,\gamma_{t}%
}(u)\mid^{2}du\mid^{\frac{1}{2}}]=0,$ we have
\[
Z^{\gamma_{t},x}(s)=\sigma(W_{s}^{\gamma_{t}},X^{\gamma_{t},x}(s),Y^{\gamma
_{t},x}(s),Z^{\gamma_{t},x}(s))\nabla_{x}u(W_{s}^{\gamma_{t}},X^{\gamma_{t}%
,x}(s))+D_{zz}u(W_{s}^{\gamma_{t}},X^{\gamma_{t},x}(s))\quad P-a.s.
\]
This completes the proof.\quad$\Box$

\section{The related path-dependent partial differential equation}

In this section, we relate FBSDE (2.1), (2.2) to the following path-dependent
partial differential equation:%
\begin{equation}%
\begin{array}
[c]{l}%
D_{t}u(\gamma_{t},x)+\mathcal{L}u(\gamma_{t},x)+tr[\nabla_{x}D_{z}u(\gamma
_{t},x)\sigma(\gamma_{t},x,u(\gamma_{t},x),v(\gamma_{t},x))]+\frac{1}%
{2}tr[D_{zz}u(\gamma_{t},x)]\\
=h(\gamma_{t},x,u(\gamma_{t},x),v(\gamma_{t},x)),\\
v(\gamma_{t},x)=\nabla_{x}u(\gamma_{t},x)\sigma(\gamma_{t},x,u(\gamma
_{t},x),v(\gamma_{t},x))+D_{z}u(\gamma_{t},x),\\
u(\gamma_{T},x)=g(\gamma_{T},x),\quad\gamma_{T}\in{\Lambda}^{n},
\end{array}
\tag{4.1}%
\end{equation}
where $\nabla$ is the gradient operator and
\[
\mathcal{L}u=\frac{1}{2}tr[(\sigma\sigma^{T})(\gamma_{t},x,u,v)\nabla
_{xx}u]+\langle b(\gamma_{t},x,u,v)\nabla_{x}u\rangle.
\]

\begin{theorem}
Suppose Assumptions $2.2,$ $3.2$ and $3.3$ hold. If $u$ belongs to
$\mathbb{C}^{1.2.2}(\Lambda\times\mathbb{R})$ and $(u,v)$ is the solution of
equation $(4.1)$ such that $(u,v)$ is uniformly Lipschitz continuous and
bounded by $C(1+|x|+\parallel\gamma_{t}\parallel)$, then we have $u(\gamma
_{t},x)=Y^{\gamma_{t},x}(t)$, for each $(\gamma_{t},x)\in\Lambda
\times\mathbb{R}$, where $(X^{\gamma_{t},x}(s),Y^{\gamma_{t},x}(s),Z^{\gamma
_{t},x}(s))_{t\leq s\leq T}$ is the unique solution of FBSDE (2.1), (2.2).
\end{theorem}

\noindent\textbf{{Proof.}} By the assumptions of this theorem, we know that
$b(\gamma_{t},x,u(\gamma_{t}.x),v(\gamma_{t},x))$ and $\sigma(\gamma
_{t},x,u(\gamma_{t},x),v(\gamma_{t},x))$ is uniformly Lipschtiz continuous.
Then the following SDE has a uniqueness solution.%
\[%
\begin{array}
[c]{rl}%
dX^{\gamma_{t},x}(s)= & b(W_{s}^{\gamma_{t}},X^{\gamma_{t},x}(s),u({W_{s}%
^{\gamma_{t}}},X^{\gamma_{t},x}(s)),v({W_{s}^{\gamma_{t}}},X^{\gamma_{t}%
,x}(s)))ds\\
& +\sigma(W_{s}^{\gamma_{t}},X^{\gamma_{t},x}(s),u({W_{s}^{\gamma_{t}}%
},X^{\gamma_{t},x}(s)),v({W_{s}^{\gamma_{t}}},X^{\gamma_{t},x}(s)))dW(s),\\
X(t)= & x,\quad s\in\lbrack t,T].
\end{array}
\]
Set
\[
(Y(s),Z(s))=(u({W_{s}^{\gamma_{t}}},X^{\gamma_{t},x}(s)),v({W_{s}^{\gamma_{t}%
}},X^{\gamma_{t},x}(s))),\quad t\leq s\leq T.
\]
We have%
\[%
\begin{array}
[c]{rl}%
dX^{\gamma_{t},x}(s)= & b(W_{s}^{\gamma_{t}},X^{\gamma_{t},x}%
(s),Y(s),Z(s))ds+\sigma(W_{s}^{\gamma_{t}},X^{\gamma_{t},x}%
(s),Y(s),Z(s))dW(s),\\
X(t)= & x,\quad s\in\lbrack t,T].
\end{array}
\]
Note that $u$ solves equation (4.1). Applying Functional $It\hat{o}$ formula
to $Y(s)=u({W_{s}^{\gamma_{t}}},X^{\gamma_{t},x}(s))$, we get%
\[%
\begin{array}
[c]{rl}%
dY(s)= & h(W_{s}^{\gamma_{t}},X^{\gamma_{t},x}(s),Y(s),Z(s))dr+Z(s)dW(s),\\
Y(T)= & g(W_{T}^{\gamma_{t}},X^{\gamma_{t},x}(T))\quad s\in\lbrack t,T].
\end{array}
\]
Thus, by the uniqueness and existence theorem of the FBSDE (Theorem 2.2), we
have the result.\quad$\Box$

Now we prove the converse to the about result.

\begin{theorem}
Under Assumptions $2.2$, $3.2$ and $3.3$, the function $u(\gamma
_{t},x)=Y^{\gamma_{t},x}(t)$ is the unique $\mathbb{C}^{1,2,2}(\Lambda
\times\mathbb{R}^{m})$-solution of the path-dependent PDE $(4.1)$.
\end{theorem}

\noindent\textbf{{Proof.}} We only study the one dimensional case. By
Corollary $3.1,$ $u\in\mathbb{C}^{0,2,2}(\Lambda)$. Let $\delta>0$ be such
that $t+\delta\leq T$. By Lemma $3.2,$ we get
\[
u(X_{t+\delta}^{\gamma_{t},x})=Y^{\gamma_{t},x}(t+\delta).
\]
Hence
\[
u(\gamma_{t,t+\delta},x)-u(\gamma_{t},x)=u(\gamma_{t,t+\delta}%
,x)-u(W_{t+\delta}^{\gamma_{t}},X^{\gamma_{t},x}(t+\delta))+u(W_{t+\delta
}^{\gamma_{t}},X^{\gamma_{t},x}(t+\delta))-u(\gamma_{t},x).
\]
Similarly as the proof of Theorem $3.3$, we obtain%
\[%
\begin{array}
[c]{rl}
& u(\gamma_{t,t+\delta},x)-u(W_{t+\delta}^{\gamma_{t}},X^{\gamma_{t}%
,x}(t+\delta))\\
= & \lim_{n\rightarrow\infty}[u^{n}(\gamma_{t,t+\delta},x)-u^{n}(W_{t+\delta
}^{\gamma_{t}},X^{\gamma_{t},x}(t+\delta))]\\
& +\int_{t}^{t+\delta}h(W_{s}^{\gamma_{t}},X^{\gamma_{t},x}(s),Y^{\gamma
_{t},x}(s),Z^{\gamma_{t},x}(s))ds+\int_{t}^{t+\delta}Z^{\gamma_{t},x}(s)dW(s).
\end{array}
\]
Following Lemma $3.1$ and Theorem $3.2$ of Pardoux and Peng \cite{Pardoux.E 2}
and Theorem $4.4$ of Peng and Wang \cite{Peng S 3}, we deduce that%
\[%
\begin{array}
[c]{rl}
& u^{n}(\gamma_{t,t+\delta},x)-u^{n}(W_{t+\delta}^{\gamma_{t}},X^{\gamma
_{t},x}(t+\delta))\\
= & \int_{t}^{t+\delta}D_{s}u^{n}({\gamma_{t,s}},x)ds-\int_{t}^{t+\delta}%
D_{s}u^{n}(W_{s}^{\gamma_{t}},X^{\gamma_{t},x}(s))ds-\int_{t}^{t+\delta}%
D_{z}u^{n}(W_{s}^{\gamma_{t}},X^{\gamma_{t},x}(s))dW(s)\\
& -\frac{1}{2}\int_{t}^{t+\delta}D_{xx}u^{n}(W_{s}^{\gamma_{t}},X^{\gamma
_{t},x}(s))ds-\frac{1}{2}\int_{t}^{t+\delta}\nabla_{xx}u^{n}(W_{s}^{\gamma
_{t}},X^{\gamma_{t},x}(s))(\sigma^{n}(W_{s}^{\gamma_{t}},X^{\gamma_{t}%
,x}(s),Y^{\gamma_{t},x}(s),Z^{\gamma_{t},x}(s)))^{2}ds\\
& -\int_{t}^{t+\delta}D_{z}\nabla_{x}u^{n}(W_{s}^{\gamma_{t}},X^{\gamma_{t}%
,x}(s))\sigma^{n}(W_{s}^{\gamma_{t}},X^{\gamma_{t},x}(s),Y^{\gamma_{t}%
,x}(s),Z^{\gamma_{t},x}(s))ds\\
& -\int_{t}^{t+\delta}\nabla_{x}u^{n}(W_{s}^{\gamma_{t}},X^{\gamma_{t}%
,x}(s))b^{n}(W_{s}^{\gamma_{t}},X^{\gamma_{t},x}(s),Y^{\gamma_{t}%
,x}(s),Z^{\gamma_{t},x}(s))ds\\
& -\int_{t}^{t+\delta}\nabla_{x}u^{n}(W_{s}^{\gamma_{t}},X^{\gamma_{t}%
,x}(s))\sigma^{n}(W_{s}^{\gamma_{t}},X^{\gamma_{t},x}(s),Y^{\gamma_{t}%
,x}(s),Z^{\gamma_{t},x}(s))dW(s).
\end{array}
\]
Thus, by the dominated convergence theorem,%
\begin{equation}%
\begin{array}
[c]{rl}
& u(\gamma_{t,t+\delta},x)-u(\gamma_{t},x)\\
= & -\int_{t}^{t+\delta}D_{z}u(W_{s}^{\gamma_{t}},X^{\gamma_{t},x}%
(s))dW(s)-\frac{1}{2}\int_{t}^{t+\delta}D_{xx}u(W_{s}^{\gamma_{t}}%
,X^{\gamma_{t},x}(s))ds\\
& -\frac{1}{2}\int_{t}^{t+\delta}\nabla_{xx}u(W_{s}^{\gamma_{t}},X^{\gamma
_{t},x}(s))(\sigma(W_{s}^{\gamma_{t}},X^{\gamma_{t},x}(s),Y^{\gamma_{t}%
,x}(s),Z^{\gamma_{t},x}(s)))^{2}ds\\
& -\int_{t}^{t+\delta}D_{z}\nabla_{x}u(W_{s}^{\gamma_{t}},X^{\gamma_{t}%
,x}(s))\sigma(W_{s}^{\gamma_{t}},X^{\gamma_{t},x}(s),Y^{\gamma_{t}%
,x}(s),Z^{\gamma_{t},x}(s))ds\\
& -\int_{t}^{t+\delta}\nabla_{x}u(W_{s}^{\gamma_{t}},X^{\gamma_{t}%
,x}(s))\sigma(W_{s}^{\gamma_{t}},X^{\gamma_{t},x}(s),Y^{\gamma_{t}%
,x}(s),Z^{\gamma_{t},x}(s))dW(s)\\
& -\int_{t}^{t+\delta}\nabla_{x}u(W_{s}^{\gamma_{t}},X^{\gamma_{t}%
,x}(s))b(W_{s}^{\gamma_{t}},X^{\gamma_{t},x}(s),Y^{\gamma_{t},x}%
(s),Z^{\gamma_{t},x}(s))ds\\
& +\int_{t}^{t+\delta}h(W_{s}^{\gamma_{t}},X^{\gamma_{t},x}(s),Y^{\gamma
_{t},x}(s),Z^{\gamma_{t},x}(s))ds+\int_{t}^{t+\delta}Z^{\gamma_{t}%
,x}(s)dW(s)+\lim_{n\rightarrow\infty}C^{n},
\end{array}
\tag{4.2}%
\end{equation}
where
\[
C^{n}=\int_{t}^{t+\delta}D_{s}u^{n}({\gamma_{t,s}},x)ds-\int_{t}^{t+\delta
}D_{s}u^{n}(W_{s}^{\gamma_{t}},X^{\gamma_{t},x}(s))ds.
\]
Since $u^{n}(\gamma_{t},x)\in\mathbb{C}_{l,lip}^{0,2,2}(\Lambda\times
\mathbb{R})$, by Lemma $3.3$ we get\newline%
\[
\mid D_{s}u^{n}(\gamma_{t,s},x)-D_{s}u^{n}(W_{s}^{\gamma_{t},x},X^{\gamma
_{t},x}(s))\mid\leq c(\parallel\gamma_{t,s}-W_{s}^{\gamma_{t},x}\parallel+\mid
X^{\gamma_{t},x}(s)-x\mid)
\]
for some constant $c$ only depending on $C,T,\gamma_{t}$ and $k$. Hence
\[
\mid C^{n}\mid\leq c\delta(\sup_{s\in{[t,t+\delta]}}\mid W^{\gamma_{t}%
}(s)-\gamma_{t}(t)\mid+\mid X^{\gamma_{t},x}(s)-x\mid).
\]
Taking expectation on both sides of $(4.2)$, we have%
\[%
\begin{array}
[c]{rl}
& \lim_{\delta\rightarrow0}\frac{u(\gamma_{t,t+\delta},x)-u(\gamma_{t}%
,x)}{\delta}\\
= & -\mathcal{L}u(\gamma_{t},x)-(\nabla_{x}{}D_{z}u(\gamma_{t},x)\sigma
(\gamma_{t},x,u(\gamma_{t},x),v(\gamma_{t},x))+\frac{1}{2}D_{xx}u(\gamma
_{t},x))+h(\gamma_{t},x,u(\gamma_{t},x),v(\gamma_{t},x)),
\end{array}
\]
where%
\[
v(\gamma_{t},x)=\nabla_{x}u(\gamma_{t},x)\sigma(\gamma_{t},x,u(\gamma
_{t},x),v(\gamma_{t},x))+D_{z}u(\gamma_{t},x).
\]
Thus, $u(\gamma_{t},x)\in\mathbb{C}^{1,2,2}(\Lambda\times\mathbb{R})$
satisfies the equation $(4.1)$.\quad$\Box$

\renewcommand{\refname}

\end{document}